\documentclass[3p]{elsarticle}
\usepackage{amssymb}
\usepackage{bm}
\usepackage{subcaption}
\usepackage{amsthm}
\usepackage{comment}
\usepackage{romannum}
\usepackage{enumitem}
\usepackage{makecell}
\usepackage[fleqn]{amsmath}

\usepackage{algorithm}
\usepackage{algpseudocode}
\usepackage{nomencl}
\makenomenclature
\usepackage{booktabs} 
\usepackage{array}
\usepackage{multirow}
\usepackage{tabularx}
\usepackage{caption}
\usepackage{siunitx}
\AtBeginDocument{\pagenumbering{arabic}}
\usepackage{url} 
\usepackage{xurl}
\biboptions{sort&compress} 
\usepackage{hyperref}
\hypersetup{pdfauthor={Yuxin Xia}}
\allowdisplaybreaks
\usepackage{xcolor}
\usepackage{longtable}
\usepackage{threeparttable}
\usepackage[utf8]{inputenc}
\usepackage[T1]{fontenc}
\hypersetup{
  colorlinks,
  citecolor=violet,
  linkcolor=red,
  urlcolor=blue}
\usepackage{graphicx} 
\usepackage{tikz}
\newcommand*\circled[1]{\tikz[baseline=(char.base)]{
            \node[shape=circle,draw,inner sep=2pt] (char) {#1};}}
\algnewcommand\algorithmicinput{\textbf{Input:}}
\algnewcommand\INPUT{\item[\algorithmicinput]}
\algdef{SE}[IF]{If}{EndIf}[1]
  {\algorithmicif\ #1\ \algorithmicthen}{\algorithmicend\ \algorithmicif}

\begin{document}

\begin{frontmatter}

\title{Lifetime Profit-Maximising Co-optimisation of Multi-Service Stacking for Battery Storage}

\author[inst1,instOMS]{Yuxin Xia \corref{cor1}}
\cortext[cor1]{Corresponding author}
\author[inst1]{Frederik Schiele}
\author[inst1]{Yihong Zhou}
\author[inst1]{Volkan Kumtepeli}
\author[inst1,instOMS]{David Howey}
\author[inst2]{Iacopo Savelli}
\author[inst1,instOMS]{Thomas Morstyn}
\affiliation[inst1]{organization={Department of Engineering Science, University of Oxford},
            city={Oxford},
            postcode={OX1 3PJ}, 
            country={UK}}
\affiliation[instOMS]{organization={Oxford Martin School, University of Oxford},
            city={Oxford},
            postcode={OX1 3BD}, 
            country={UK}}
\affiliation[inst2]{organization={Centre for Research on Geography, Resources, Environment, Energy $\&$ Networks, Bocconi University},
            city={Milano},
            postcode={20136}, 
            country={Italy}}

\begin{abstract}
Grid-scale battery energy storage can generate revenue by stacking services across electricity and frequency response markets, yet identifying the lifetime profit-maximising stacking strategy remains challenging. Decisions across services are coupled through shared battery system capacity, constrained by system operator energy management rules, and further shaped by product-specific technical requirements that govern system operation, degradation, and lifetime profitability.
This paper presents an ageing-aware receding-horizon framework for co-optimising multi-service stacking that explicitly captures product-specific characteristics and state-of-energy compliance rules. The framework is applied to the Great Britain market, where storage operators can stack electricity trading with multiple dynamic frequency response services procured through the newly introduced `Enduring Auction Capability’ platform under energy management requirements imposed by the National Energy System Operator.
Using real market data, we demonstrate that degradation modelling, discount rate, and battery ageing jointly govern both lifetime value and optimal stacking strategy. Accounting for ageing increases lifetime revenue by up to 32\% relative to a degradation-agnostic benchmark, while higher-fidelity ageing modelling can provide a further revenue improvement of up to 16\% over simpler formulations. Lower discount rates favour strategies that balance calendar and cycling ageing, while higher discount rates favour aggressive operation. Across all strategies, services responding to positive frequency deviations are consistently preferred over those responding to negative deviations.
Lifetime profit is maximised by adapting the optimal service mix as the battery ages. In the early years, the optimal strategy prioritises baseline discharge and fast-acting frequency response, while service provision diversifies as degradation increases.
\end{abstract}
\end{frontmatter}

\section{Introduction}\label{section1}
The United Kingdom has reinforced its statutory commitment to achieving a net-zero economy by 2050 through the 2019 amendment to the Climate Change Act \cite{UK_Climate_Change_Act}, complemented by the government’s Clean Power 2030 Action Plan \cite{DESNZ_2024_CleanPower2030}. The resulting growth in renewable generation, together with declining system inertia, has increased the operational challenge of maintaining system balance. In response to this, battery energy storage systems (BESSs) have emerged as a key source of power system flexibility and reliability \cite{HUANG2023100126}. They are highly controllable and capable of rapid response, making them well suited to providing balancing and frequency response services. Declining production costs have driven a rapid expansion of installed grid-scale battery capacity in Great Britain (GB) in recent years \cite{statkraft2025roleBESS,NREL_2023}, with BESSs increasingly relied upon by the National Energy System Operator (NESO) to support system reliability and security of supply \cite{NESO_FES2025}.

To address rising frequency response requirements under low-inertia system conditions, NESO introduced a new suite of dynamic frequency response (DFR) services from 2020 onwards, procured through the new Enduring Auction Capability platform \cite{NESO2025DynamicServices,NESO_EAC}. These services comprise Dynamic Regulation (DR), Dynamic Moderation (DM), and Dynamic Containment (DC), each with distinct technical characteristics. They have since fully replaced legacy frequency response arrangements and now represent an important revenue source for grid-scale BESS participation in GB's ancillary service markets \cite{NESO2025DynamicServices}. 

Across international markets, frequency response procurement arrangements and service designs differ. PJM in the United States adopts a multi-product co-optimisation framework that enables the concurrent clearing of multiple regulation services \cite{PJM_2018_Regulation}, whereas Nordic markets rely on separate, service-specific procurement mechanisms \cite{ENTSOE_2024_Regulation}. In the Australian NEM, energy and frequency control ancillary services are co-optimised in dispatch, while individual frequency control services are enabled and priced separately in each region \cite{AEMO2025FCASGuide}. ERCOT, the Texas electricity market, has only recently transitioned towards real-time co-optimisation under its 2025 RTC+B reform \cite{AlKez2023}. In contrast, the GB Enduring Auction Capability platform adopts a unified auction framework with co-optimisation across multiple frequency service categories and endogenous capacity allocation. The DR, DM, DC, and reserve products are jointly cleared within a single optimisation, supported by advanced market features such as order stacking and splitting, negative-price acceptance, and algorithmic over-procurement if it reduces total clearing cost \cite{NESO_EAC_ppt}. In its inaugural year, this integrated design delivered substantial efficiency gains: total ancillary service costs fell by 32\%, while procured volumes rose by 52\%, and the weighted average price declined from £3.99/MWh to £1.77/MWh \cite{NSIDE2025PoweringGBGrid}. These outcomes suggest that the Enduring Auction Capability platform design offers useful insights for the reform of frequency regulation markets in other countries.

While this integrated market design improves system-level efficiency, it also imposes more stringent operational requirements on energy-limited technologies such as BESSs. Dynamic frequency response services are cleared day-ahead with availability payments for real-time activation, and sustained service delivery over the contract period induces cumulative changes in the state-of-energy (SoE), creating an inherent operational constraint. Without appropriate management, these changes can progressively drive the asset into operating conditions where the available energy is insufficient to deliver the full contracted response \cite{NESO_SOE_monitor,NGESO_SOE_ppt,WU2024114413}. To address this risk, as part of the June 2024 Article 18 consultation on dynamic response services under the Electricity Balancing Regulation, NESO introduced new SoE management and monitoring rules \cite{NESO_SOE_monitor}, with penalty mechanisms for non-compliance enforced from April 2025 \cite{neso2025_announce}. If the required SoE thresholds set by NESO are not met, providers must restore their energy position through baseline submissions that reserve power capacity for energy recovery \cite{NESO_SOE_monitor,UKbattery_SoE_mana}. In this way, DFR services can continue to be delivered in a sustained and compliant manner over the contract period.

Alongside DFR services, storage assets can provide a range of market and system services through revenue stacking, including imbalance settlement, reserve provision, energy arbitrage, and peak shaving \cite{NESO2025DynamicServices,MOY2021100065,11050105,SEWARD2022108292}. In particular, arbitrage enables BESS to exploit temporal price differentials by charging during low-price periods and discharging during high-price periods. Given prior evidence that stacking multiple services can yield substantially higher revenues for grid-scale BESS \cite{ENGLBERGER2020100238}, combinations that integrate energy market participation with frequency response services, such as arbitrage and DFR, are of particular interest. 

However, arbitrage and DFR provision are intrinsically coupled \cite{10478754}. This coupling arises because the upward and downward flexibility delivered in DFR services is defined relative to the battery’s baseline energy level, which is also determined by arbitrage decisions in electricity markets. As an energy-limited asset, the feasible provision of both services is therefore jointly constrained by the battery’s energy and power capacity. This interdependence is further reinforced by NESO SoE management requirements imposed on DFR providers, which require maintaining SoE within prescribed bounds and recovering energy following service delivery to ensure contract compliance across trading intervals \cite{NESO_SOE_monitor}. This requires submitting baseline charging or discharging schedules in electricity markets, directly linking DFR provision to arbitrage decisions. As a result, arbitrage and DFR decisions are strongly interdependent and cannot be optimised independently.

Driven by recent cost reductions and technological advances, lithium-ion batteries, particularly lithium iron phosphate chemistries, are becoming the dominant technology in GB. Projections indicate that lithium-ion batteries will account for around 90\% of GB grid-scale BESS deployments by 2030 \cite{Faraday_RhoMotion_2023_UK_BESS}. Within lithium-ion technologies, the share of lithium iron phosphate chemistries is expected to increase from approximately 60\% in 2022 to over 70\% by 2030, driven by favourable cost, performance, and degradation characteristics \cite{Faraday_RhoMotion_2023_UK_BESS}. Nevertheless, lithium-ion batteries are inherently subject to degradation driven by multiple internal ageing mechanisms, including solid electrolyte interphase growth, loss of active material, and loss of lithium inventory, leading to capacity and power fade \cite{BIRKL2017373,HAN2019100005,VETTER2005269}. In particular, capacity fade reduces the usable energy of the system and constrains the lifetime economic performance of BESS assets \cite{Kumtepeli_Howey_2022}, whereas power fade generally does not constitute a binding constraint owing to the moderate C-rates typical of stationary applications. In practice, lithium-ion batteries are commonly regarded as having reached their economic end-of-life (EoL) when the remaining usable capacity, quantified by the state-of-health (SoH), falls below a threshold typically set at 80\% of the initial rated capacity in line with manufacturer warranty specifications \cite{7488267}. This threshold is widely associated with the onset of accelerated, non-linear degradation beyond the `knee point', accompanied by increased operational uncertainty and elevated safety risks \cite{HAN2019100005,Attia2022Knees}.

From an operational perspective, the stress factors driving battery degradation are commonly classified into calendar ageing and cycle ageing. Calendar ageing is governed primarily by time, SoC, and temperature, while cycle ageing is driven by usage-related factors such as charge and discharge rates, depth of discharge, full cycle equivalents, cutoff voltages, and temperature \cite{COLLATH2022105634,MPC_battery,en12060999}. While some ageing drivers are intrinsic to cell design and manufacturing quality, many operational stress factors are directly controllable through dispatch and energy management decisions. This enables the development of ageing-aware operating strategies that mitigate degradation and improve the lifetime economic performance of grid-scale BESSs \cite{COLLATH2022105634}.

For a grid-scale BESS asset, operators must determine how battery capacity should be allocated across multiple revenue streams over time. This requires an optimal multi-service stacking strategy across different electricity market products, such as energy arbitrage, frequency response, reserve, and balancing services, to maximise combined market revenues while accounting for service-specific requirements and battery degradation \cite{10886993}. The preferred strategy may also change as the battery ages, since each service has different revenue potential, technical requirements, utilisation patterns, and degradation impacts. Therefore, revenue-stacking strategies should be evaluated not only by their short-term market returns, but also by their effects on degradation, usable lifetime, and long-term asset value.

Prior studies have demonstrated that ageing-aware BESS operational optimisation can improve lifetime profitability, with existing approaches differing in the representation of battery ageing within the optimisation framework \cite{COLLATH2022105634,MPC_battery,Volkan2019}. Strategies that neglect ageing optimise short-term revenues without explicitly accounting for degradation, which can lead to operating patterns that undermine long-term performance and asset value. Rule-based strategies introduce ageing considerations implicitly via predefined operational constraints, including C-rate limits, daily cycle caps, and constrained voltage operating windows. In contrast, optimisation-based strategies represent ageing explicitly by modelling battery degradation and typically incorporate it as a penalty term within the objective function \cite{MPC_battery}. This formulation provides a systematic way of addressing the trade-off between revenue and degradation cost.

The operational optimisation of grid-scale BESS depends on both short-term economic returns and the long-term impacts of degradation on usable lifetime. Although a growing body of literature examines ageing-aware BESS operation (see Table~\ref{tab:LR}), most studies either consider only cycle ageing or restrict the analysis to a limited set of applications, typically focusing on a single service such as arbitrage and neglecting the optimisation of multi-service stacking strategies \cite{MPC_battery,9250459,WANKMULLER201756,RENIERS201891,EJEH20221957}. A representative study is presented in Collath et al.\ \cite{MPC_battery}, where a semi-empirical degradation model for commercial lithium iron phosphate cells is linearised using experimental data \cite{cyc_ageing_LFP} and embedded within a mixed-integer linear programming formulation. Although this approach achieved a 30\% improvement in lifetime profitability for arbitrage in the German intra-day market compared with an ageing-unaware strategy, the analysis was limited to a single service. This restriction is significant given prior evidence that stacking multiple services can yield substantially higher revenues for grid-scale BESS \cite{ENGLBERGER2020100238}. Existing studies that consider multiple applications either neglect calendar degradation or lack validation against high-fidelity battery models, which remains essential for the deployment of ageing-aware strategies in real-world BESS operation \cite{ENGLBERGER2020100238,UKbattery_SoE_mana,MARIAUD2017466}.

\begin{table}[t]
\centering
\caption{Summary of the literature review.}
\label{tab:LR}
\renewcommand{\arraystretch}{1.3}
\setlength{\tabcolsep}{4pt}
\footnotesize

\begin{tabular}{%
p{1.4cm}
p{1.8cm}
p{3.0cm}
p{1.5cm}
p{1.6cm}
p{1.8cm}
p{1.4cm}
p{1.5cm}}
\toprule

\makecell[t]{\textbf{Ref.}} &
\makecell[t]{\textbf{Market}\\\textbf{scope}} &
\makecell[t]{\textbf{Ageing-aware}\\\textbf{approach}} &
\makecell[t]{\textbf{Cycle}\\\textbf{ageing}} &
\makecell[t]{\textbf{Calendar}\\\textbf{ageing}} &
\raisebox{0.35ex}{\makecell[t]{\textbf{NESO SoE}\\\textbf{management}}} &
\makecell[t]{\textbf{MPC}\\\textbf{approach}} &
\makecell[t]{\textbf{Region}} \\
\midrule

\cite{ENGLBERGER2020100238} &
EA, FR, PS &
OB: Cycle ageing cost &
Yes &
No &
N/A &
Partially &
Germany \\

\cite{MPC_battery} &
EA &
OB: Cycle \& calendar ageing cost &
Yes &
Yes &
N/A &
Yes &
Germany \\

\cite{9250459} &
EA &
OB: Cycle \& calendar ageing cost &
Yes &
Yes &
N/A &
Yes &
Germany \\

\cite{WANKMULLER201756} &
EA &
OB: Cycle ageing cost, RB: SoC range limit &
Yes &
No &
N/A &
No &
United States (MISO) \\

\cite{10045057} &
FR &
OB: Cycle ageing cost &
Yes &
No &
N/A &
No &
United States (PJM) \\

\cite{RENIERS201891} &
EA &
OB: Cycle \& calendar ageing cost &
Yes &
Yes &
N/A &
No &
Belgium \\

\cite{EJEH20221957} &
EA &
Not considered &
No &
No &
Not considered &
No &
GB \\

\cite{MARIAUD2017466} &
EA, FR &
Not considered &
No &
No &
Not considered &
No &
GB \\

\cite{MARTINS2021111938} &
EA, RP, FR &
OB: Cycle ageing cost &
Yes &
No &
Not considered &
No &
GB \\

\cite{WU2024114413} &
DFR (DC service only) &
Not considered &
No &
No &
Simplified &
No &
GB \\

\cite{UKbattery_SoE_mana} &
EA, DFR, IMB &
RB: Daily cycle limit &
Yes &
No &
Simplified &
No &
GB \\

This paper &
EA, DFR &
OB: Cycle \& calendar ageing cost, RB: Daily cycle limit &
Yes &
Yes &
Detailed &
Yes &
GB \\

\bottomrule
\end{tabular}

\vspace{2mm}
\begin{minipage}{\textwidth}
\footnotesize
\raggedright
\textit{Abbreviations:} EA, energy arbitrage; FR, frequency response; 
RP, reserve provision; IMB, imbalance settlement; DFR, dynamic frequency response; 
PS, peak shaving; RB, rule-based; OB, optimisation-based; GB, Great Britain.
\end{minipage}
\end{table}

In the UK context, the introduction of DFR services under stringent NESO SoE requirements introduces additional operational complexity that remains insufficiently addressed in existing ageing-aware BESS optimisation studies. As summarised in Table~\ref{tab:LR}, only two GB-focused studies \cite{WU2024114413,UKbattery_SoE_mana} explicitly consider DFR participation, and both remain limited in scope. The first study \cite{UKbattery_SoE_mana} proposes a three-stage optimisation strategy for BESS participation across energy arbitrage, DFR and imbalance settlement markets. However, it restricts each time block to a single DFR service and therefore does not capture multi-service stacking across the DC, DM and DR products. In addition, battery degradation is represented through a rule-based daily cycle cap rather than an explicit cost-based formulation, and no validation against a detailed battery model is provided. The second study \cite{WU2024114413} focuses on market clearing for a virtual BESS providing DC services under a simplified representation of NESO SoE management. While this provides useful insights into frequency response provision, it considers only a single DFR service and omits other revenue streams associated with DR and DM services and arbitrage opportunities in the electricity market. Although SoE management is incorporated in these studies, it is represented in a simplified manner that does not capture the full operational implications of NESO requirements, particularly the need for baseline energy adjustments and inter-temporal coupling between service provision and energy scheduling. As a result, these approaches may overestimate achievable revenues and misrepresent feasible operational strategies. More importantly, existing GB-focused studies do not establish how the optimal allocation of EA and multiple DFR services should evolve over the battery lifetime as service characteristics, operating constraints, and battery condition interact. They also do not capture how cumulative degradation changes long-term operational behaviour and economic performance.


To address these gaps, we develop an ageing-aware receding-horizon co-optimisation framework for grid-scale BESSs that jointly stacks electricity arbitrage and multiple frequency response services under the new Enduring Auction Capability platform. The framework provides decision support for battery operators and asset managers by identifying which services should be prioritised at different stages of battery life to maximise long-term asset value while managing degradation. The main contributions of this study are summarised as follows:
\begin{enumerate}
    \item A co-optimisation framework for grid-scale BESS that captures ageing-aware multi-service stacking across two electricity market products and six DFR services, explicitly representing service-specific technical characteristics aligned with NESO design, cross-service coupling through shared baseline schedules, battery energy and power constraints, and detailed NESO SoE management requirements.

    \item A systematic quantification of how ageing-aware operational strategies and discount rates influence projected lifetime profitability and long-term decision-making for grid-scale BESS under GB market conditions.

    \item A lifetime perspective on ageing-aware service stacking, showing how the economically optimal allocation of multiple arbitrage and DFR services with different technical characteristics shifts over the battery lifetime to maximise lifetime profitability under cumulative degradation.
\end{enumerate}


\section{Dynamic Frequency Response Services and Energy Management Framework in GB}
\label{sec:background}
This section outlines the regulatory and operational background informing the optimisation framework developed in this study. Section~\ref{sec:bg:DFR} reviews the DFR services and their associated technical specifications. Section~\ref{sec:bg:soe} then introduces the SoE management rules implemented by NESO, which govern the admissible energy trajectories of BESSs participating in the DFR market. The service names and abbreviations used in this paper follow the official terminology defined in the NESO Service Terms \cite{NESO_Service_Terms}.

\subsection{Dynamic frequency response services}\label{sec:bg:DFR}
In November 2023, GB restructured its procurement framework for DFR services with the introduction of the Enduring Auction Capability platform to maintain system frequency within the band of $50 \pm 0.5$~Hz \cite{NESO_EAC,neso_eac_auction_results}. Under this framework, NESO procures six DFR products grouped within DR, DM, and DC, each with low- and high-frequency response variants, namely `DRL'/`DRH', `DML'/`DMH', and `DCL'/`DCH'. Dynamic regulation provides continuous but slower frequency correction, while dynamic moderation supports system stability during periods of heightened volatility. Dynamic containment complements these pre-fault services by delivering fast response to large system disturbances~\cite{NESO2025DynamicServices}. The technical specifications and procurement volumes in 2024 for each service are summarised in Table~\ref{tab:dfr_technical_requirements}.

Dynamic frequency response services are procured in alignment with the six standard electricity forward agreement blocks, each spanning four hours and commencing at fixed times of 23:00, 03:00, 07:00, 11:00, 15:00, and 19:00 (GMT), with each block comprising eight 30-minute settlement periods \cite{NESO_MARKET_DESIGN_REPORT}. Market participants may bid into and deliver multiple DFR services simultaneously, including both low-frequency and high-frequency services, which can be contracted on either a single-directional or bi-directional basis. Accepted capacity is remunerated through \textit{availability} payments settled at the corresponding clearing price in \pounds/MW/h. In addition, each contract specifies a contracted response energy volume that defines the maximum energy volume that a response unit may be required to deliver under an activation instruction during the contracted service period \cite{NESO_Service_Terms}. This volume is calculated as the product of the contracted power (MW) and the associated delivery duration (h) for each electricity forward agreement block. It therefore represents the contractual maximum energy obligation per activation, and serves as a reference for compliance monitoring and SoE management by NESO.

\begin{table}[t]
\centering
\caption{Technical specifications and annual procurement data for the three DFR services in 2024~\cite{neso_eac_auction_results}.}
\renewcommand{\arraystretch}{1.3}
\setlength{\tabcolsep}{4pt}
\footnotesize
\begin{tabular}{p{1.5cm} p{2.8cm} p{2.9cm} p{2.3cm} p{3.5cm}}
\toprule
\midrule
\textbf{Dynamic \newline Service} & 
\textbf{Initiation Time} & 
\textbf{Time to Full \newline  Delivery} & 
\textbf{Delivery \newline  Duration} & 
\textbf{Procured Capacity} \\
\midrule
Description & 
Max. time between frequency change and start of response &
Max. time between \newline frequency change and full response delivery &
Time of sustained \newline response delivery &
Total response capacity \newline procured by NESO in 2024 \\
\midrule
DC & 0.5\,s & 1\,s & 15\,min & 
H: 2,712,488 MW\\
& & & & 
L: 2,453,760 MW\\
\addlinespace[2pt]
DM & 0.5\,s & 1\,s & 30\,min & 
H: 409,248 MW\\
& & & & 
L: 348,281 MW \\
\addlinespace[2pt]
DR & 2\,s & 10\,s & 60\,min & 
H: 757,299 MW\\
& & & & 
L: 660,028 MW \\
\midrule
\bottomrule
\end{tabular}
\label{tab:dfr_technical_requirements}
\end{table}

\begin{figure}[t] 
  \centering
  \includegraphics[width=0.8\columnwidth]{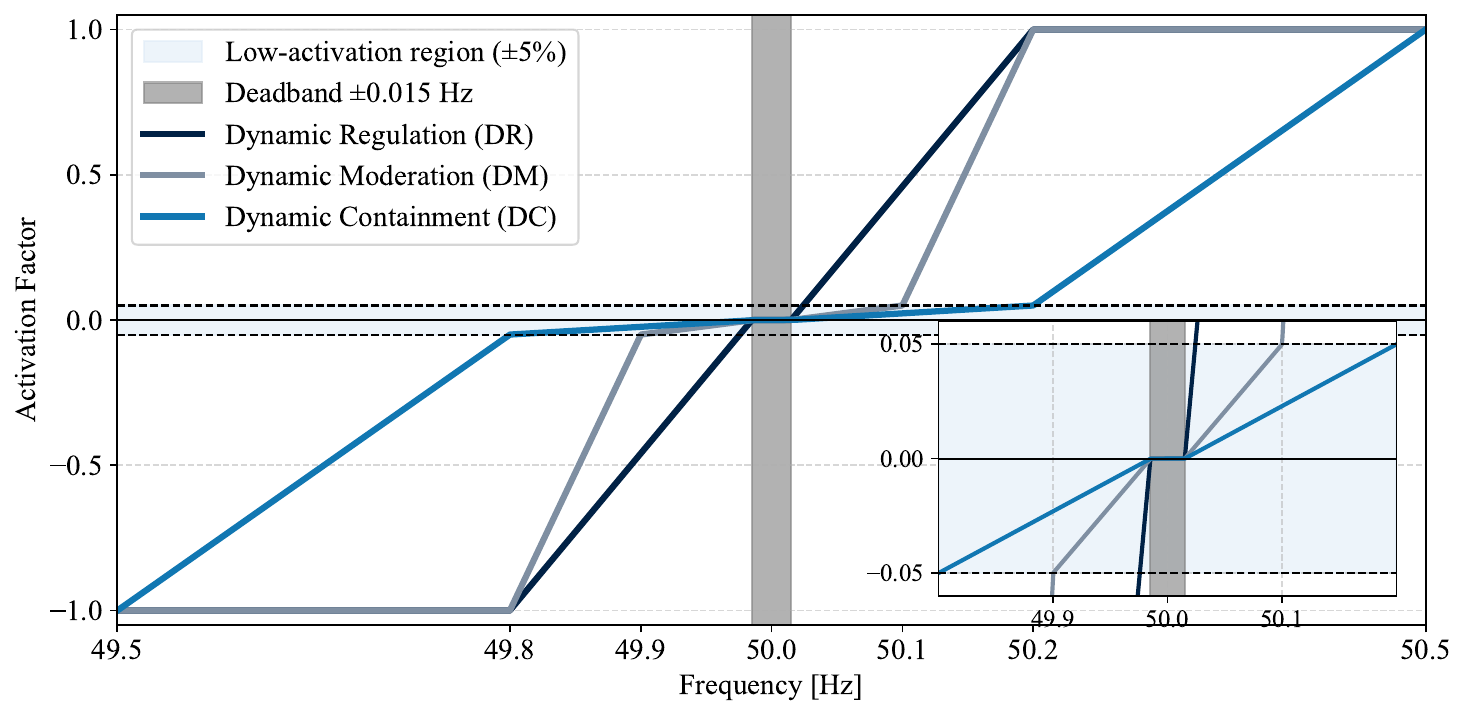}
  \caption{Illustration of the activation factors $\alpha_t$ for DC, DM and DR. The activation factor describes the share of the contracted power that has to be delivered at a certain grid frequency. DFR high-frequency services are activated at frequencies higher than 50 Hz (right side of the graph) while low-frequency services are activated during times of frequencies below 50 Hz (left side of the graph). Activation factors for DFR low-frequency services are negative, indicating power absorption by the BESS, while activation factors for DFR high-frequency services are positive, indicating power injection into the BESS.}\label{fig:fre_act}
\end{figure}

During service delivery, the activated proportion of contracted response power is determined by real-time frequency deviation according to service-specific activation profiles. As shown in Fig.~\ref{fig:fre_act}, each DFR service $s \in \mathcal{S} = \{\mathrm{DCH, DCL, DMH, DML, DRH, DRL}\}$ is described by a piecewise-linear activation function $g^{s}(f_t)$, mapping system frequency $f_t$ to an activation factor $\alpha_t^{s}$. A deadband of ±0.015 Hz around 50 Hz (49.985–50.015 Hz) applies within which no activation occurs. Outside this range, activation increases according to service-specific thresholds and ramp rates. High frequency services are triggered when $f_t>50.015$ Hz, leading to $\alpha_t^{s}>0$ and battery charging. Conversely, low-frequency services activate when $f_t>49.985$ Hz, yielding $\alpha_t^{s}<0$ and battery discharging. For example, at 49.9 Hz, a unit contracted for 4 MW of DRL with $\alpha_t^{\mathrm{DRL}} = -0.5$ delivers 2 MW discharge. Over a 30-minute settlement period, this corresponds to 1 MWh of frequency response energy.

\subsection{NESO Energy management rules} \label{sec:bg:soe}
For energy-limited storage assets providing DFR services under availability-based remuneration, NESO imposes SoE management rules to ensure sustained delivery of the contracted response throughout the contracted service period, i.e., the time window during which the provider is contractually required to deliver the service (typically a 4-hour trading block) \cite{NESO_SOE_monitor,NGESO_SOE_ppt}. For such units, realised frequency response activation progressively changes the SoE, which may constrain a unit's ability to maintain its committed response power over time. To preserve continued deliverability, NESO implements an energy recovery framework that regulates post-activation SoE restoration. For units holding bi-directional contracts, additional power reserve requirements are imposed at the day-ahead bidding stage to ensure sufficient operational flexibility for recovery.

The energy recovery framework is implemented through a dynamic minimum SoE requirement that serves as the minimum compliance threshold monitored by NESO. The minimum SoE requirement is initially set equal to the contracted response energy volume for each 4-hour trading block, representing the energy capability that must be preserved for the contracted service period. It is subsequently updated at the start of each settlement period to reflect cumulative service delivery and recovery actions. Realised frequency response energy delivered during the preceding settlement period decreases the minimum SoE requirement, whereas scheduled energy recovery adjustment volume increases it. 

At the beginning of each settlement period, the reported SoE is assessed against the prevailing minimum SoE requirement using performance monitoring data, as part of NESO’s behavioural and compliance checks designed to ensure proper energy management and service integrity. When the SoE falls below the minimum SoE requirement, providers must submit baseline charge or discharge decisions in the electricity market at the earliest opportunity. Such updates are submitted in the next settlement period and take effect only after gate closure, typically one hour (two settlement periods) later. Consequently, recovery of energy delivered in a given settlement period can only occur four settlement periods later. For example, energy delivered in settlement period 1 can only be restored from settlement period 5 onward, since the updated baseline for settlement period 5 must be submitted in settlement period 2. This operational constraint is referred to hereafter as the 4-settlement-period delay. Through this mechanism, the minimum SoE requirement defines a dynamic lower bound on SoE, embedding the energy recovery requirement into the compliance constraint and ensuring sustained deliverability of the contracted response.

We present two examples of this in \ref{app:neso_soe} (Fig.~\ref{fig:neso_case1} and Fig.~\ref{fig:neso_case2}), including detailed explanations of the resulting allowable SoE range under NESO rules and the interaction between the minimum SoE requirement and recovery volumes. The mathematical formulation of the SoE management framework is provided in Section~\ref{sec:mpc:opt:soe_mana}.

\section{Proposed model predictive control framework}\label{sec:mpc}
\begin{figure}[t] 
  \centering
  \includegraphics[width=0.9\columnwidth]{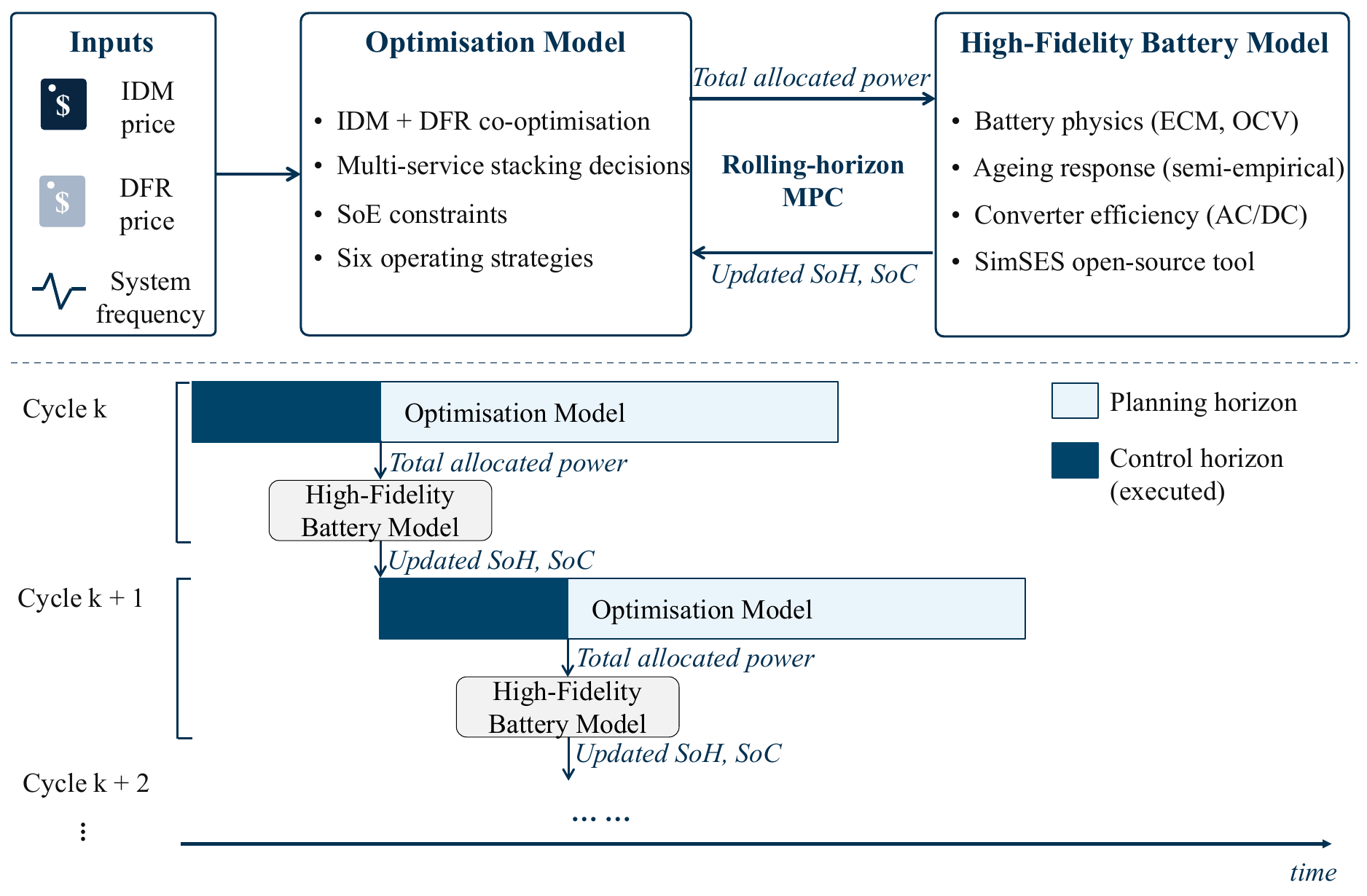}
  \caption{Rolling-horizon MPC framework for ageing-aware multi-service BESS scheduling. The optimisation model determines market-service allocations and passes the resulting total power trajectory to a high-fidelity battery model for validation. The updated SoH and SoC are then fed back to initialise the next optimisation cycle.}
\label{fig:mpc_frame}
\end{figure}

Fig.~\ref{fig:mpc_frame} illustrates the overall model predictive control framework for ageing-aware battery operation, and Fig.~\ref{fig:mpc_frame} shows the rolling-horizon principle. It comprises two key components; an optimisation model and a high-fidelity battery model. This structure enables systematic benchmarking of different ageing-aware optimisation formulations on a digital twin of the BESS and closely reflects a realistic implementation. The optimisation model determines the operating strategy and the high-fidelity battery model feeds back essential operational states, including SoC and SoH. Consequently, the framework provides a controlled environment in which operational strategies can be designed and validated prior to deployment on a real BESS. 


The optimisation model determines the optimal capacity allocation and energy schedule for EA and DFR participation within each rolling horizon. Inputs include GB electricity prices, DFR clearing prices, and system frequency data, together with simulated internal battery state information provided by the high-fidelity battery model, including SoC and SoH. Using these inputs, the optimisation model determines the operating decisions that maximise market profit under a simplified BESS representation. To reflect the different revenue strategies that may be adopted by BESS operators and asset managers in practice, six optimisation strategies are proposed. These strategies are implemented through variants of the same rolling-horizon optimisation framework, with different objective-function terms and ageing-cost representations. They range from short-term revenue maximisation without explicit ageing considerations to ageing-aware operation that accounts for degradation costs and long-term battery value. At each optimisation step, a two-day planning horizon is solved within the rolling-horizon framework illustrated in Fig.~\ref{fig:mpc_frame}. Only the first day of the optimal dispatch, corresponding to the control horizon, is implemented as the AC power profile for daily operation and passed to the high-fidelity battery model for validation and closed-loop assessment.

Based on the optimal dispatch sequence determined by the optimisation model, the high-fidelity model simulates the real-world response of the BESS. For this purpose we employed the open-source SimSES modelling framework for stationary BESS \cite{SimSES_paper}. It includes an equivalent-circuit representation of the cells together with open-circuit voltage curves, a semi-empirical ageing model to simulate capacity fade and resistance growth, and power converter models that capture AC–DC efficiency and associated losses. SimSES simulates a realistic BESS response to validate the optimisation results and update the corresponding SoC and SoH. These results are then fed back to the optimisation model, allowing it to initialise the subsequent optimisation step. This process was repeated iteratively until the assumed or actual end of life. %
In what follows, the optimisation problem is presented in Section~\ref{sec:mpc:opt}, with the high-fidelity battery model described in Section~\ref{sec:mpc:hfbm}.

\subsection{Optimisation models} \label{sec:mpc:opt}
In this study, six optimisation model configurations are evaluated within a co-optimised multi-service framework, capturing participation in electricity markets through baseline charging and discharging for arbitrage, and in DFR services including DC, DM, and DR across both high- and low-frequency products, under NESO SoE management constraints. Table~\ref{tab:optim_models} summarises these configurations by listing the six ageing models alongside their corresponding objective functions and constraint sets. The remainder of this section introduces each ageing model and the SoE management rule in detail.

\begin{table}[t]
\centering
\caption{Overview of the six optimisation models, including objective functions and constraint sets.}
\begin{tabular}{l c c c c}
\toprule
\midrule
\textbf{Model} & 
\textbf{Objective} & 
\textbf{\shortstack{Ageing-aware \\ Constraints}} & 
\textbf{\shortstack{NESO SoE \\ Constraints}} & 
\textbf{\shortstack{Other \\ Constraints}} \\
\midrule
\textit{No ageing} & \eqref{eq:obj_no_ageing} & / & \eqref{eq:neso_soe_mana} & \eqref{eq:base}  \\
\textit{Cycle limit} & \eqref{eq:obj_no_ageing} & \eqref{eq:cyc_limit} & \eqref{eq:neso_soe_mana} & \eqref{eq:base} \\
\textit{L cyc.} & \eqref{eq:obj_ageing_cyc} & \eqref{eq:cycle_ageing_L} & \eqref{eq:neso_soe_mana} & \eqref{eq:base} \\
\textit{L cal. \& cyc.} & \eqref{eq:obj_ageing_cal_cyc} & \eqref{eq:calendar_ageing_L}, \eqref{eq:cycle_ageing_L} & \eqref{eq:neso_soe_mana} & \eqref{eq:base} \\
\textit{PL cyc.} & \eqref{eq:obj_ageing_cyc} & \eqref{eq:cycle_ageing_PL} & \eqref{eq:neso_soe_mana} & \eqref{eq:base} \\
\textit{PL cal. \& cyc.} & \eqref{eq:obj_ageing_cal_cyc} & \eqref{eq:calendar_ageing_PL}, \eqref{eq:cycle_ageing_PL} & \eqref{eq:neso_soe_mana} & \eqref{eq:base} \\
\midrule
\bottomrule
\end{tabular}
\label{tab:optim_models}
\end{table}

\subsubsection{No ageing}
The \textit{no ageing} model is a formulation that ignores battery degradation. Neither the objective function nor the constraints contain terms or limits related to degradation. The optimisation focuses purely on short-term revenue maximisation subject only to market signals and technical feasibility. 
\begin{equation}
\max \quad 
\Pi^{\mathrm{E}}
+  \Pi^{\mathrm{DFR}}
= \sum_{t\in \mathcal{T}_{\mathrm{SP}}} \big(p^{\mathrm{base,dis}}_{t}-p^{\mathrm{base,ch}}_{t}\big)\pi^{\mathrm{IDM}}_{t} \Lambda^{\mathrm{SP}}+ \sum_{t\in \mathcal{T}}\sum_{s\in\mathcal{S}}
p^{s}_{t} \pi^{s}_{t}\Lambda^{\mathrm{EFA}}\label{eq:obj_no_ageing}
\end{equation}
The objective function~\eqref{eq:obj_no_ageing} maximises the total revenue obtained from the electricity market, \( \Pi^{\mathrm{E}} \), and from the DFR market, \( \Pi^{\mathrm{DFR}} \). The first component captures the profit arising from the baseline charging and discharging schedule across all settlement periods within the optimisation horizon $\mathcal{T}_{\mathrm{SP}}$. This revenue is computed by multiplying the net baseline power with the intra-day electricity market price \( \pi^{\mathrm{IDM}}_{t} \) published by EPEX SPOT~\cite{epexspot_gbdata}, and the settlement period duration \( \Lambda^{\mathrm{SP}} = 0.5\,\mathrm{h} \). The second component captures the revenue associated with providing DFR. For each optimisation timestep \( t \in \mathcal{T} \) with resolution \( \Delta t \), the contracted power \( p_t^s \) is multiplied by the service clearing price \( \pi_t^s \) published by NESO~\cite{neso_eac_auction_results} and the EFA duration \( \Lambda^{\mathrm{EFA}} = 4\,\mathrm{h} \) to compute the revenue from the DFR market.

Eqs.\ \eqref{eq:base} specify the technical constraints of the BESS operation:
\begin{subequations}\label{eq:base}
\begin{align}
&p_t^{\mathrm{tot}} = p_t^{\mathrm{base,ch}} - p_t^{\mathrm{base,dis}}
+\big( \alpha_t^{\mathrm{DCH}} p_t^{\mathrm{DCH}} + \alpha_t^{\mathrm{DMH}}p_t^{\mathrm{DMH}} + \alpha_t^{\mathrm{DRH}}p_t^{\mathrm{DRH}}\big) \notag\\
&\qquad\quad- \big(\alpha_t^{\mathrm{DCL}}p_t^{\mathrm{DCL}} + \alpha_t^{\mathrm{DML}}p_t^{\mathrm{DML}} + \alpha_t^{\mathrm{DRL}}p_t^{\mathrm{DRL}}\big),
&& \forall t\in\mathcal{T}, \label{eq:base:power_balance}\\
&0 \leq p_t^{s} \leq P^{\max},
&& \forall t\in\mathcal{T}_\mathrm{EFA},s\in\mathcal{S}, \label{eq:base:service_cap}\\
&z_t^{\mathrm{base,ch}} + z_t^{\mathrm{base,dis}} \leq 1,
&& \forall t\in\mathcal{T}_\mathrm{SP}, \label{eq:base:baseline_excl}\\
&0 \leq p_t^{\mathrm{base,ch}} \leq z_t^{\mathrm{base,ch}} P^{\max},
&& \forall t\in\mathcal{T}_\mathrm{SP}, \label{eq:base:baseline_ch}\\
&0 \leq p_t^{\mathrm{base,dis}} \leq z_t^{\mathrm{base,dis}} P^{\max},
&& \forall t\in\mathcal{T}_\mathrm{SP}, \label{eq:base:baseline_dis}\\
&p_t^{\mathrm{tot}} = p_t^{\mathrm{tot,ch}} - p_t^{\mathrm{tot,dis}},
&& \forall t\in\mathcal{T}, \label{eq:base:split_tot}\\
&z_t^{\mathrm{tot,ch}} + z_t^{\mathrm{tot,dis}} \leq 1,
&& \forall t\in\mathcal{T}, 
\label{eq:base:total_excl}\\
&0 \leq p_t^{\mathrm{tot,ch}} \le z_t^{\mathrm{tot,ch}}P^{\max},
&& \forall t\in\mathcal{T}, \label{eq:base:ch_box}\\
&0 \leq p_t^{\mathrm{tot,dis}} \le z_t^{\mathrm{tot,dis}}P^{\max},
&& \forall t\in\mathcal{T}, \label{eq:base:dis_box}\\
&\mathrm{SoC}_{t} = \mathrm{SoC}_{t-1}
+ \frac{\Delta t}{E^{\mathrm{batt}}}\!\left(\eta\, p_{t-1}^{\mathrm{tot,ch}}
- \frac{1}{\eta}\, p_{t-1}^{\mathrm{tot,dis}}\right),
&& \forall t\in\mathcal{T}\setminus\{0\}, \label{eq:base:soc_dyn}\\
&\mathrm{SoC}_{0} = \mathrm{SoC}^{\mathrm{ini}}, \label{eq:base:soc_init}\\
&0 \leq \mathrm{SoC}_{t} \le 1,
&& \forall t\in\mathcal{T}. \label{eq:base:soc_box}
\end{align}
\end{subequations}
Constraint~\eqref{eq:base:power_balance} defines the total power exchanged with the grid at each optimisation timestep $\Delta t$ as the sum of the baseline and DFR service schedules. The activation factor \( \alpha_t^s \in [0,1] \) represents the proportion of capacity activated in response to system frequency deviations. This factor can be derived from system frequency data~\cite{neso_system_frequency}. Constraint~\eqref{eq:base:service_cap} ensures that the power allocated to each DFR service remains within the BESS’s maximum power limit \( P^{\max} \). Constraints~\eqref{eq:base:baseline_excl}--\eqref{eq:base:baseline_dis} enforce mutual exclusivity and power bounds on the baseline charging and discharging actions. 
Constraints~\eqref{eq:base:split_tot}--\eqref{eq:base:dis_box} impose the same exclusivity and power limits on the total BESS power, which aggregates baseline EA and DFR services. Constraints~\eqref{eq:base:soc_dyn}--\eqref{eq:base:soc_box} define the evolution of the battery's SoC over time. The SoC is updated recursively based on the net charging and discharging energy. Constraint~\eqref{eq:base:soc_dyn} tracks the normalised stored energy (SoE) of the BESS. The initial SoC is set to $\mathrm{SoC}^{\mathrm{ini}}$ as in~\eqref{eq:base:soc_init}, and the SoC is constrained between 0 and 1 by~\eqref{eq:base:soc_box}. In battery scheduling models for power-system applications, the battery state is commonly represented as normalised stored energy, as in \eqref{eq:base:soc_dyn}, although this quantity is often referred to as SoC. In electrochemical modelling, by contrast, SoC is defined on a charge basis.  In this paper, the normalised SoE from the scheduling model is used as a proxy for the SoC input required by the ageing model. This assumption is reasonable for the operating conditions considered here, where average C-rates are moderate and, for LFP batteries, the open-circuit voltage is relatively flat over much of the practical operating SoC range \cite{BRAUN2026239471}.

\subsubsection{Cycle limit}
The \textit{cycle limit} model is formulated by imposing a daily full cycle equivalent limit, denoted by $n^{\mathrm{cyc,max}}$. Constraint~\eqref{eq:cyc_limit} enforces this daily full cycle equivalent cap within the rolling horizon. The cycle count is computed as the sum of charging and discharging half-cycles.
\begin{align}\label{eq:cyc_limit}
n_t^{\mathrm{cyc}} =
\sum_{\tau \in \mathcal{H}^{\mathrm{daily}}_{\Delta t}(t)}
\left(
\frac{p_\tau^{\mathrm{tot,ch}} + p_\tau^{\mathrm{tot,dis}}}{2E^{\mathrm{batt}}}
\cdot \frac{\Delta t}{3600}
\right)
\le n^{\mathrm{cyc,max}},\quad \forall\, t \in \mathcal{T}_{\mathrm{daily}}.
\end{align}

\subsubsection{Cycle and calendar ageing}
\begin{figure}[t] 
  \centering
  \includegraphics[width=\columnwidth]{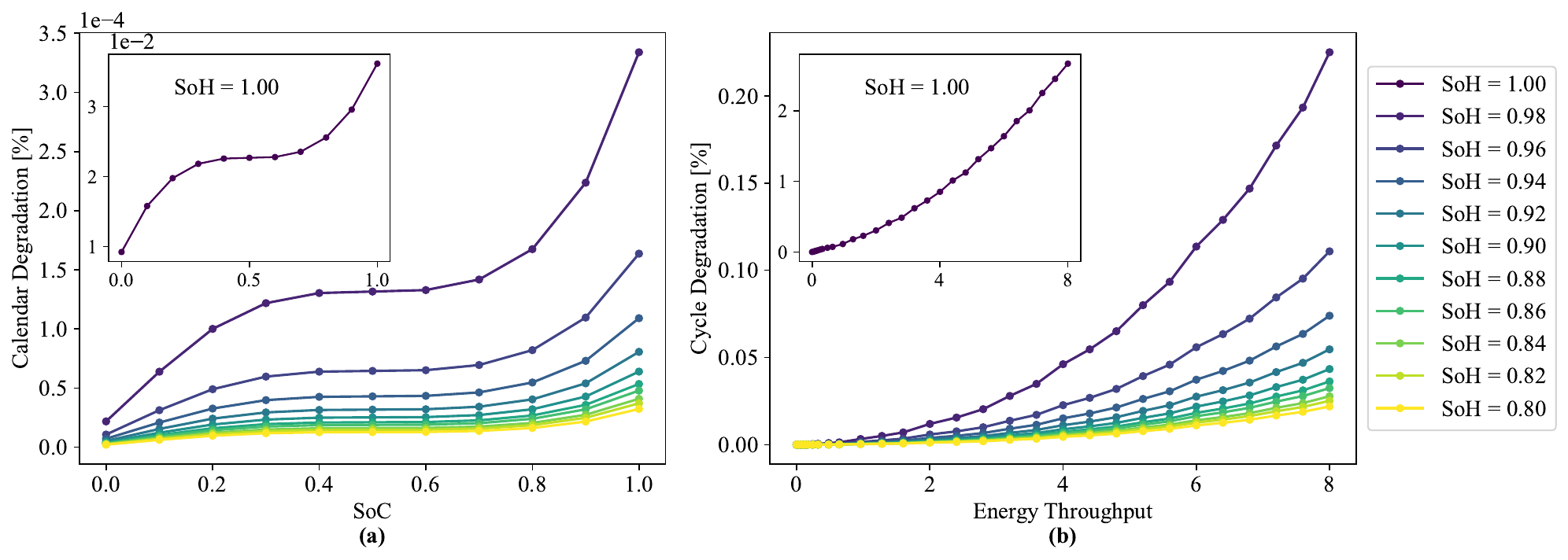}
  \caption{Piecewise-linear approximations of the semi-empirical ageing model for an LFP battery: (a) calendar ageing over a 15-min period as a function of SoC, and (b) cycle ageing over a 4-hour period as a function of energy throughput. Each curve corresponds to a different SoH level, representing a different level of accumulated capacity loss. The SoH = 1.00 case is shown in the inset to avoid it dominating the vertical scale. Further details of the underlying ageing model are provided in \cite{MPC_battery,cyc_ageing_LFP,calen_ageing_LFP}.}
  \label{fig:ageing_plot}
\end{figure}

In addition to the cycle-limit constraint, we incorporate a semi-empirical degradation cost model that accounts for both calendar and cycle ageing 
\cite{cyc_ageing_LFP,calen_ageing_LFP}, together with a piecewise linear ageing cost model \cite{MPC_battery}. Both calendar and cycle degradation exhibit a square-root dependence on past capacity loss. In the optimisation model, this non-linear relationship is approximated using piecewise linear representations \cite{MPC_battery}, as illustrated in Fig.~\ref{fig:ageing_plot}. This degradation model is incorporated into the objective function as a cost term that reflects the economic impact of projected capacity loss. We evaluated four ageing cost models across two modelling dimensions: formulation complexity, comparing piecewise-linear and single-segment affine approximations, and ageing scope, comparing cycle-only ageing with combined calendar and cycle ageing.

For the cycle-only formulation, the objective function is
\begin{equation}
\max \quad
\Pi^{\mathrm{E}}
+\Pi^{\mathrm{DFR}}
-C^{\mathrm{cyc}},
\label{eq:obj_ageing_cyc}
\end{equation}
where the cycle ageing cost is given by
\begin{equation}
C^{\mathrm{cyc}}
=
\frac{E^{\mathrm{batt}} \cdot \pi^{\mathrm{loss}}}
{1-\mathrm{SoH}^{\mathrm{EoL}}}
\left(
\sum_{t\in\mathcal{T}_{\mathrm{4h}}}
q_t^{\mathrm{loss,cyc}}
\right).
\label{eq:cost_ageing_cyc}
\end{equation}

When both calendar and cycle ageing are considered, the objective function includes separate cost terms for the two degradation mechanisms:
\begin{equation}
\max \quad
\Pi^{\mathrm{E}}
+\Pi^{\mathrm{DFR}}
-C^{\mathrm{cyc}}
-C^{\mathrm{cal}},
\label{eq:obj_ageing_cal_cyc}
\end{equation}
where the calendar ageing cost is given by
\begin{equation}
C^{\mathrm{cal}}
=
\frac{E^{\mathrm{batt}} \cdot \pi^{\mathrm{loss}}}
{1-\mathrm{SoH}^{\mathrm{EoL}}}
\left(
\sum_{t\in\mathcal{T}_{\mathrm{15m}}}
q_t^{\mathrm{loss,cal}}
\right),
\label{eq:cost_ageing_cal}
\end{equation}
and \(C^{\mathrm{cyc}}\) is defined in~\eqref{eq:cost_ageing_cyc}. Here, \(E^{\mathrm{batt}}\) denotes the battery energy capacity, \(\pi^{\mathrm{loss}}\) represents the assumed economic value of capacity loss, and \(\mathrm{SoH}^{\mathrm{EoL}}\) denotes the end-of-life SoH threshold. The degradation terms \(q_t^{\mathrm{loss,cyc}}\) and \(q_t^{\mathrm{loss,cal}}\) represent the estimated fractional capacity losses due to cycle and calendar ageing, respectively.

\subsubsection*{Piecewise Linear (PL) Calendar Ageing}
\begin{subequations}\label{eq:calendar_ageing_PL}
\begin{align}
&\overline{\mathrm{SoC}}_{t}
= \sum_{\tau \in \mathcal{H}^\mathrm{15m}_{\Delta t}(t)} \mathrm{SoC}_{\tau} \frac{\Delta t}{900},
&&\forall t \in \mathcal{T}_{15\mathrm{m}}, \label{eq:cal_PL_1}\\
&\lambda^{\text{cal}}_{t,j} \leq 1,  && \forall t \in \mathcal{T}_\mathrm{15m},  j \in \mathcal{J}, \label{eq:cal_PL_2} \\[6pt]
&\sum_{j \in \mathcal{J}} \lambda^{\text{cal}}_{t,j} = 1, && \forall t \in \mathcal{T}_\mathrm{15m}, \label{eq:cal_PL_3} \\
&\sum_{j \in \mathcal{J}} \lambda^{\text{cal}}_{t,j} \cdot X^{\text{cal}}_j = \overline{\mathrm{SoC}}_{t},  && \forall t \in \mathcal{T}_\mathrm{15m}, \label{eq:cal_PL_4} \\
&\sum_{j \in \mathcal{J}} \lambda^{\text{cal}}_{t,j} \cdot Z^{\text{cal}}_j = q_{t}^{\text{loss,cal}},  && \forall t \in \mathcal{T}_\mathrm{15m}. \label{eq:cal_PL_5}
\end{align}
\end{subequations}

\subsubsection*{Linear (L) Calendar Ageing}
\begin{subequations}
\label{eq:calendar_ageing_L} 
\begin{align} 
&\overline{\mathrm{SoC}}_{t}
= \sum_{\tau \in \mathcal{H}^\mathrm{15m}_{\Delta t}(t)} \mathrm{SoC}_{\tau} \frac{\Delta t}{900},
&&\forall t \in \mathcal{T}_{\mathrm{15m}}, \label{eq:cal_L_1}\\
&q^{\text{loss,cal}}_{t} =a_{\text{cal}}\:\overline{\mathrm{SoC}}_{t}+b_{\text{cal}}, &&\forall t\in\mathcal{T}_{\mathrm{15m}}. \label{eq:cal_L_2}
\end{align} 
\end{subequations}

The piecewise linear formulation of calendar ageing follows~\cite{MPC_battery} and is defined by~\eqref{eq:calendar_ageing_PL}. Constraint~\eqref{eq:cal_PL_1} computes the 15-minute average state of charge \( \overline{\text{SoC}}_t \), which serves as the input to the calendar ageing function shown in Fig.~\ref{fig:ageing_plot}(a). This function is defined over a set of breakpoints \( \{X^{\text{cal}}_j, Z^{\text{cal}}_j\}_{j \in \mathcal{J}} \), where \( X^{\text{cal}}_j \) denotes discrete SoC levels and \( Z^{\text{cal}}_j \) the corresponding calendar degradation values. The interpolation is implemented via the SOS2 variables \( \lambda^{\text{cal}}_{t,j} \) in~\eqref{eq:cal_PL_2}--\eqref{eq:cal_PL_5}, which act as convex combination weights. In particular, \eqref{eq:cal_PL_3}--\eqref{eq:cal_PL_4} enforce that \( \overline{\text{SoC}}_t \) is represented as a weighted combination of the breakpoints, while \eqref{eq:cal_PL_5} maps these weights to the corresponding degradation value $q_{t}^{\text{loss,cal}}$. The SOS2 structure ensures that at most two adjacent \( \lambda^{\text{cal}}_{t,j} \) are non-zero, thereby restricting the solution to lie on a single segment of the piecewise linear approximation and enabling linear interpolation between neighbouring breakpoints. For comparison, a simplified linearised (L) calendar ageing model is introduced in~\eqref{eq:calendar_ageing_L}, where the degradation is approximated as a single-segment affine function of the average SoC. This formulation removes the need for SOS2 constraints, reducing computational complexity at the expense of modelling accuracy.

\subsubsection*{Piecewise Linear (PL) Cycle Ageing}
\begin{subequations}
\label{eq:cycle_ageing_PL}
\begin{align}
&e^{\mathrm{ch}}_{t} = \sum_{\tau \in \mathcal{H}^\mathrm{4h}_{\Delta t}(t)} p_{\tau}^{\text{tot,ch}}  \frac{\Delta t}{3600}, && \forall t \in \mathcal{T}_\mathrm{4h}, \label{eq:cyc_PL_1} \\
&e^{\mathrm{dis}}_{t} = \sum_{\tau \in \mathcal{H}^\mathrm{4h}_{\Delta t}(t)} p_{\tau}^{\text{tot,dis}}  \frac{\Delta t}{3600}, && \forall t \in \mathcal{T}_\mathrm{4h}, \label{eq:cyc_PL_2} \\
&\lambda^{\text{cyc}}_{t,k} \leq 1, && \forall t \in \mathcal{T}_\mathrm{4h}, \forall k \in \mathcal{K}, \label{eq:cyc_PL_3} \\
&\sum_{k \in \mathcal{K}} \lambda^{\text{cyc,ch}}_{t,k} = 1, && \forall t \in \mathcal{T}_\mathrm{4h}, \label{eq:cyc_PL_4} \\
&\sum_{k \in \mathcal{K}} \lambda^{\text{cyc,dis}}_{t,k} = 1, && \forall t \in \mathcal{T}_\mathrm{4h}, \label{eq:cyc_PL_5} \\
&\sum_{k \in \mathcal{K}} \lambda^{\text{cyc,ch}}_{t,k} \cdot X^{\text{cyc}}_k = e^{\mathrm{ch}}_{t}, && \forall t \in \mathcal{T}_\mathrm{4h}, \label{eq:cyc_PL_6} \\
&\sum_{k \in \mathcal{K}} \lambda^{\text{cyc,dis}}_{t,k} \cdot X^{\text{cyc}}_k = e^{\mathrm{dis}}_{t}, && \forall t \in \mathcal{T}_\mathrm{4h}, \label{eq:cyc_PL_7} \\
&\sum_{k \in \mathcal{K}} \lambda^{\text{cyc,ch}}_{t,k} \cdot Z^{\text{cyc}}_k = q_{t}^{\text{loss,cyc,ch}}, && \forall t \in \mathcal{T}_\mathrm{4h}, \label{eq:cyc_PL_8} \\
&\sum_{k \in \mathcal{K}} \lambda^{\text{cyc,dis}}_{t,k} \cdot Z^{\text{cyc}}_k = q_{t}^{\text{loss,cyc,dis}}, && \forall t \in \mathcal{T}_\mathrm{4h}, \label{eq:cyc_PL_9} \\
&q^{\text{loss,cyc}}_{t} = q^{\text{loss,cyc,ch}}_{t} + q^{\text{loss,cyc,dis}}_{t}, && \forall t \in \mathcal{T}_\mathrm{4h}. \label{eq:cyc_PL_10}
\end{align}
\end{subequations}

\subsubsection*{Linear (L) Cycle Ageing}
\begin{subequations}\label{eq:cycle_ageing_L} \begin{align} 
&e^{\mathrm{ch}}_{t} = \sum_{\tau \in \mathcal{H}^\mathrm{4h}_{\Delta t}(t)} p_{\tau}^{\text{tot,ch}}  \frac{\Delta t}{3600}, && \forall t \in \mathcal{T}_\mathrm{4h}, \label{eq:cyc_L_1} \\
&e^{\mathrm{dis}}_{t} = \sum_{\tau \in \mathcal{H}^\mathrm{4h}_{\Delta t}(t)} p_{\tau}^{\text{tot,dis}}  \frac{\Delta t}{3600}, && \forall t \in \mathcal{T}_\mathrm{4h}, \label{eq:cyc_L_2} \\
&q^{\text{loss,cyc,ch}}_{t}=a^{\text{cyc,ch}}\,e^{\mathrm{ch}}_{t}+b^{\text{cyc,ch}}, &&\forall t\in\mathcal{T}_\mathrm{4h},\label{eq:cyc_L_3} \\
&q^{\text{loss,cyc,dis}}_{t}=a^{\text{cyc,dis}}\,e^{\mathrm{dis}}_{t}+b^{\text{cyc,dis}}, &&\forall t\in\mathcal{T}_\mathrm{4h},\label{eq:cyc_L_4} \\
&q^{\text{loss,cyc}}_{t} =q^{\text{loss,cyc,ch}}_{t}+q^{\text{loss,cyc,dis}}_{t}, &&\forall t\in\mathcal{T}_\mathrm{4h}. \label{eq:cyc_L_5} 
\end{align} \end{subequations}

The PL cycle ageing model in~\eqref{eq:cycle_ageing_PL} follows the same SOS2-based piecewise linear interpolation principle as the PL calendar ageing formulation. Instead of the average SoE, the inputs to the piecewise linear function are the total charging and discharging energy over each 4-hour interval, \( e^{\mathrm{ch}}_{t} \) and \( e^{\mathrm{dis}}_{t} \), as defined in~\eqref{eq:cyc_PL_1}--\eqref{eq:cyc_PL_2}, which serve as inputs to the cycle ageing function shown in Fig.~\ref{fig:ageing_plot}(b). The degradation is evaluated using SOS2 variables \( \lambda^{\text{cyc,ch}}_{t,k} \) and \( \lambda^{\text{cyc,dis}}_{t,k} \), which act as convex combination weights over the breakpoint set \( \{X^{\text{cyc}}_k, Z^{\text{cyc}}_k\}_{k \in \mathcal{K}} \). Here, \( X^{\text{cyc}}_k \) represents discretised energy throughput levels and \( Z^{\text{cyc}}_k \) the corresponding cycle degradation values. Through~\eqref{eq:cyc_PL_3}--\eqref{eq:cyc_PL_9}, the energy throughput is mapped onto the piecewise linear approximation, and the associated degradation is obtained via linear interpolation between two adjacent breakpoints. The resulting charging and discharging degradation components are then combined to yield the total cycle degradation $q^{\text{loss,cyc}}_{t}$ in~\eqref{eq:cyc_PL_10}. A simplified linear approximation of this relationship is provided in~\eqref{eq:cycle_ageing_L}.

\subsubsection{System operator energy management rules} \label{sec:mpc:opt:soe_mana}

\begin{table}[t]
\centering
\footnotesize
\caption{Key abbreviations used in the NESO SoE management constraints.}
\label{tab:soe_management_abbreviations}
\begin{tabularx}{\textwidth}{p{1.8cm}p{3.5cm}X}
\toprule
\textbf{Abbreviation} & \textbf{Full name} & \textbf{Meaning and explanation} \\
\midrule

EFA & Electricity Forward Agreement &
A 4-hour delivery block used for DFR service procurement. Each block spans four hours and commences at fixed times of 23:00, 03:00, 07:00, 11:00, 15:00, and
19:00 (GMT). \\

SP & Settlement Period &
A 30-minute settlement interval within an electricity forward agreement block. Each electricity forward agreement block contains eight settlement periods. \\

CREV & Contracted Response Energy Volume &
The contracted response energy volume that the unit must be capable of delivering at the start of a contracted electricity forward agreement period. It is calculated from the contracted quantity and the delivery duration of the relevant DFR service.\\

ER & Energy Recovery &
The maximum energy volume in MWh that the unit is able to recover within one settlement period after service activation. It is calculated as 20\% of the contracted response energy volume. \\

ERAV & Energy Recovery Adjustment Volume &
The recovery adjustment volume applied to the relevant minimum SoE requirement after service activation. \\

FRE & Frequency Response Energy &
The realised energy delivered due to frequency response activation during a settlement period. \\

MG & Margin &
The available SoE margin relative to the relevant contracted response energy volume. This captures how much of the recovery requirement can be absorbed while still satisfying the contracted response energy volume, rather than restoring a unit to its pre-event SoE level. \\

MSER & Minimum State-of-Energy Requirement &
The minimum SoE level that a unit must satisfy at the start of each settlement period. At the start of a 4-hour electricity forward agreement trading block, this is equal to the contracted response energy volume. It is then reduced by frequency response energy delivered through service activation and increased by the energy recovery adjustment volume, but cannot exceed the contracted response energy volume. \\

RER & Required Energy Recovery &
The recovery volume required in the current settlement period for each service direction. This is calculated as the realised frequency response energy plus any leftover recovery from the previous SP, capped by the energy recovery limit. It is then compared with the available margin to determine how much recovery can be absorbed and how much needs to be carried forward as an adjustment. \\

ABS & Absorption &
The portion of the required energy recovery covered by the available SoE margin. \\

adj0 & Scheduled adjustment volume &
The recovery adjustment volume scheduled in the current settlement period. This does not affect the SoE requirement immediately, but is carried forward over the 4-settlement-period implementation delay.\\

adj4 & Implemented adjustment volume &
The recovery adjustment volume implemented after the 4-settlement-period delay. This is the delayed effect of a previously scheduled adjustment and is applied to the minimum SoE requirement. \\

\bottomrule
\end{tabularx}
\end{table}

As discussed in Section~\ref{sec:bg:soe}, NESO applies a dynamic SoE management framework. This framework ensures that energy-limited providers of DFR services maintain sufficient usable energy to continuously deliver their contracted response volumes~\cite{NESO_SOE_monitor,NGESO_SOE_ppt,NESO_Service_Terms,NGESO_DC_Guidance}. The full set of NESO SoE management rules is formalised in~\eqref{eq:neso_soe_mana}. To keep the SoE management constraints compact, Table~\ref{tab:soe_management_abbreviations} defines the abbreviations used in the notation, including variable names, superscripts, and subscripts.
\begin{subequations}\label{eq:neso_soe_mana}
\begin{align}
& p_t^{\mathrm{DCL}} + p_t^{\mathrm{DML}} + p_t^{\mathrm{DRL}}  \leq P^{\max}  b^{\mathrm{L}}_t, && \forall t\in\mathcal{T}, \label{eq:neso_soe_mana:bL_def}\\
& p_t^{\mathrm{DCH}} + p_t^{\mathrm{DMH}} + p_t^{\mathrm{DRH}}  \leq P^{\max}  b^{\mathrm{H}}_t, && \forall t\in\mathcal{T}, \label{eq:neso_soe_mana:bH_def} \\
&r_t^{\mathrm{L}} \geq 0.2\left( \frac{15\cdot2}{60} p_t^{\mathrm{DCL}} +  \frac{30\cdot2}{60}p_t^{\mathrm{DML}} +  \frac{60\cdot2}{60}p_t^{\mathrm{DRL}}\right)b^{\mathrm{H}}_t  b^{\mathrm{L}}_t , && \forall t \in \mathcal{T}_\mathrm{EFA}, \label{eq:neso_soe_mana:reserve_L}\\
&r_t^{\mathrm{H}} \geq  0.2\left(\frac{15\cdot2}{60} p_t^{\mathrm{DCH}} + \frac{30\cdot2}{60} p_t^{\mathrm{DMH}} +  \frac{60\cdot2}{60}p_t^{\mathrm{DRH}}\right)b^{\mathrm{H}}_t  b^{\mathrm{L}}_t , && \forall t \in \mathcal{T}_\mathrm{EFA},\label{eq:neso_soe_mana:reserve_H}\\
&0 \leq r_t^{\mathrm{L}} \leq P^{\max}  b^{\mathrm{H}}_t, && \forall t \in \mathcal{T}_\mathrm{EFA},\label{eq:neso_soe_mana:reserve_L_bound1}\\
&0\leq r_t^{\mathrm{H}} \leq P^{\max}  b^{\mathrm{L}}_t, && \forall t \in \mathcal{T}_\mathrm{EFA},\label{eq:neso_soe_mana:reserve_H_bound1}\\
&r_t^{\mathrm{L}} \geq b^{\mathrm{H}}_t  (1-b^{\mathrm{L}}_t)  P^{\max}, && \forall t \in \mathcal{T}_\mathrm{EFA},\label{eq:neso_soe_mana:reserve_L_bound2}\\
&r_t^{\mathrm{H}} \geq (1-b^{\mathrm{H}}_t)  b^{\mathrm{L}}_t  P^{\max}, && \forall t \in \mathcal{T}_\mathrm{EFA},\label{eq:neso_soe_mana:reserve_H_bound2}\\
&p_t^{\mathrm{base,dis}} + p_t^{\mathrm{DCL}} + p_t^{\mathrm{DML}} + p_t^{\mathrm{DRL}} + r_t^{\mathrm{L}} \leq P^{\max} , && \forall t\in\mathcal{T}, \label{eq:neso_soe_mana:L_pmax}\\
&p_t^{\mathrm{base,ch}} + p_t^{\mathrm{DCH}} + p_t^{\mathrm{DMH}} + p_t^{\mathrm{DRH}} + r_t^{\mathrm{H}} \leq P^{\max} , && \forall t\in\mathcal{T},  \label{eq:neso_soe_mana:H_pmax} \\
&\mathrm{CREV}_t^{\mathrm{L}} = \tfrac14 p_t^{\mathrm{DCL}} + \tfrac12 p_t^{\mathrm{DML}} + p_t^{\mathrm{DRL}} , && \forall t\in \mathcal{T}_\mathrm{EFA}, \label{eq:neso_soe_mana:CREV_L}\\
&\mathrm{CREV}_t^{\mathrm{H}} = \tfrac14 p_t^{\mathrm{DCH}} + \tfrac12 p_t^{\mathrm{DMH}} + p_t^{\mathrm{DRH}} , && \forall t\in \mathcal{T}_\mathrm{EFA}, \label{eq:neso_soe_mana:CREV_H}\\
&\mathrm{ER}_t^{\mathrm{L}} = 0.2 \mathrm{CREV}_t^{\mathrm{L}},\qquad \mathrm{ER}_t^{\mathrm{H}} = 0.2 \mathrm{CREV}_t^{\mathrm{H}}, && \forall t\in \mathcal{T}_\mathrm{EFA}, \label{eq:neso_soe_mana:ER}\\
&\mathrm{FRE}_t^{\mathrm{L}}
=
\sum_{\tau \in \mathcal{H}_{\mathrm{\Delta t}}^{\mathrm{SP}}(t)}
\left(
\alpha_{\tau}^{\mathrm{DCL}} p_{\tau}^{\mathrm{DCL}}
+
\alpha_{\tau}^{\mathrm{DML}} p_{\tau}^{\mathrm{DML}}
+
\alpha_{\tau}^{\mathrm{DRL}} p_{\tau}^{\mathrm{DRL}}
\right)
\frac{\Delta t}{3600},
&&
\forall t \in \mathcal{T}_{\mathrm{SP}},\label{eq:neso_soe_mana:FRE_L_drf}\\
& \mathrm{FRE}_t^{\mathrm{H}}
=
\sum_{\tau \in \mathcal{H}_{\mathrm{\Delta t}}^{\mathrm{SP}}(t)}
\left(
\alpha_{\tau}^{\mathrm{DCH}} p_{\tau}^{\mathrm{DCH}}
+
\alpha_{\tau}^{\mathrm{DMH}} p_{\tau}^{\mathrm{DMH}}
+
\alpha_{\tau}^{\mathrm{DRH}} p_{\tau}^{\mathrm{DRH}}
\right)
\frac{\Delta t}{3600},
&& \forall\, t \in \mathcal{T}_{\mathrm{SP}},\label{eq:neso_soe_mana:FRE_H_drf}\\
&\mathrm{FRE}_{t}^{\mathrm{L}} \leq \tfrac14 p_t^{\mathrm{DCL}} + \tfrac12 p_t^{\mathrm{DML}} + \tfrac12 p_t^{\mathrm{DRL}} , && \forall t\in \mathcal{T}_\mathrm{SP}, \label{eq:neso_soe_mana:REV_L_bound}\\
&\mathrm{FRE}_{t}^{\mathrm{H}} \leq \tfrac14 p_t^{\mathrm{DCH}} + \tfrac12 p_t^{\mathrm{DMH}} + \tfrac12 p_t^{\mathrm{DRH}} , && \forall t\in \mathcal{T}_\mathrm{SP}, \label{eq:neso_soe_mana:REV_H_bound}\\
&\sum_{\tau\in\mathcal{H}_{\mathrm{SP}}^{\mathrm{EFA}}(t)} p_\tau^{\mathrm{base,dis}} \geq \delta\sum_{\tau\in\mathcal{H}_{\mathrm{SP}}^{\mathrm{EFA}}(t)}  \mathrm{FRE}_\tau^{\mathrm{H}}, && \forall t\in\mathcal{T}_\mathrm{EFA}, \label{eq:neso_soe_mana:offset_H}\\
&\sum_{\tau\in\mathcal{H}_{\mathrm{SP}}^{\mathrm{EFA}}(t)} p_\tau^{\mathrm{base,ch}} \geq \delta \sum_{\tau\in\mathcal{H}_{\mathrm{SP}}^{\mathrm{EFA}}(t)} \mathrm{FRE}_\tau^{\mathrm{L}}, && \forall t\in\mathcal{T}_\mathrm{EFA},  \label{eq:neso_soe_mana:offset_L}\\
& \mathrm{MSER}_{t}^{\mathrm{L}} = \mathrm{CREV}_{t}^{\mathrm{L}}, \qquad
  \mathrm{MSER}_{t}^{\mathrm{H}} = \mathrm{CREV}_{t}^{\mathrm{H}},
&& \forall\, t \in \mathcal{T}_{\mathrm{EFA}}^{\mathrm{SP1}}, \label{eq:neso_soe_mana:MSER_START}\\
& \mathrm{MSER}_{t}^{\mathrm{L}} \leq \mathrm{CREV}_t^{\mathrm{L}},\qquad \mathrm{MSER}_{t}^{\mathrm{H}} \leq \mathrm{CREV}_t^{\mathrm{H}}, 
&& \forall t\in \mathcal{T}_\mathrm{SP}, \label{eq:neso_soe_mana:MSER_allSP}\\
& \mathrm{MSER}_{t+1}^{\mathrm{L}} = \mathrm{MSER}_{t}^{\mathrm{L}} + \mathrm{adj4}_{t}^{\mathrm{L}} - \mathrm{FRE}_{t}^{\mathrm{L}}, 
&& \forall\, t \in \mathcal{T}_{\mathrm{EFA}}^{\mathrm{SP2-8}}, \label{eq:neso_soe_mana:MSER_cal_L}\\[4pt]
& \mathrm{MSER}_{t+1}^{\mathrm{H}} = \mathrm{MSER}_{t}^{\mathrm{H}} + \mathrm{adj4}_{t}^{\mathrm{H}} - \mathrm{FRE}_{t}^{\mathrm{H}}, 
&& \forall\, t \in \mathcal{T}_{\mathrm{EFA}}^{\mathrm{SP2-8}},\label{eq:neso_soe_mana:MSER_cal_H}\\
&E^{\mathrm{batt}} \mathrm{SoC}_t \geq \mathrm{MSER}_t^{\mathrm{L}}  , && \forall t\in \mathcal{T}_\mathrm{SP},\label{eq:neso_soe_mana:soc_L_bound}\\
& E^{\mathrm{batt}} \mathrm{SoC}_t \leq E^{\mathrm{batt}}-\mathrm{MSER}_t^{\mathrm{H}} , && \forall t\in \mathcal{T}_\mathrm{SP}, \label{eq:neso_soe_mana:soc_H_bound}\\
& \mathrm{MG}_t^{\mathrm{L}} = \max\!\left\{0, E^{\mathrm{batt}} \mathrm{SoC}_t - \mathrm{CREV}^{\mathrm{L}}_t \right\}, && \forall t\in \mathcal{T}_\mathrm{SP},\label{eq:neso_soe_mana:mg_cal_L}\\
& \mathrm{MG}_t^{\mathrm{H}} = \max\!\left\{0, E^{\mathrm{batt}} - E^{\mathrm{batt}} \mathrm{SoC}_t - \mathrm{CREV}^{\mathrm{H}}_t \right\}, && \forall t\in \mathcal{T}_\mathrm{SP},\label{eq:neso_soe_mana:mg_cal_H}\\
&\mathrm{RER}_t^{\mathrm{L}} = \min\{ \mathrm{FRE}_{t}^{\mathrm{L}} + \mathrm{left}_t^{\mathrm{L}}, \mathrm{ER}_t^{\mathrm{L}} \}, && \forall t\in \mathcal{T}_\mathrm{SP},\label{eq:neso_soe_mana:RER_cal_L}\\
&\mathrm{RER}_t^{\mathrm{H}} = \min\{\mathrm{FRE}_{t}^{\mathrm{H}} + \mathrm{left}_t^{\mathrm{H}}, \mathrm{ER}_t^{\mathrm{H}} \}, && \forall t\in \mathcal{T}_\mathrm{SP},\label{eq:neso_soe_mana:RER_cal_H}\\
& \mathrm{ABS}_t^{\mathrm{L}} = \min\left\{ \mathrm{RER}_t^{\mathrm{L}}, \mathrm{MG}_t^{\mathrm{L}}\right\}, && \forall t\in \mathcal{T}_\mathrm{SP},\label{eq:neso_soe_mana:ABSR_cal_L}\\
& \mathrm{ABS}_t^{\mathrm{H}} = \min\left\{ \mathrm{RER}_t^{\mathrm{H}}, \mathrm{MG}_t^{\mathrm{H}}\right\}, && \forall t\in \mathcal{T}_\mathrm{SP},\label{eq:neso_soe_mana:ABSR_cal_H}\\
& \mathrm{adj0}^{\mathrm{L}}_t = \mathrm{RER}_t^{\mathrm{L}} - \mathrm{ABS}_t^{\mathrm{L}},\qquad \mathrm{adj0}^{\mathrm{H}}_t = \mathrm{RER}_t^{\mathrm{H}} - \mathrm{ABS}_t^{\mathrm{H}},&& \forall t\in \mathcal{T}_\mathrm{SP},\label{eq:neso_soe_mana:ad0_def}\\
& \mathrm{adj4}_{t}^{\mathrm{L}} =
\begin{cases}
0, & \text{if } t \leq \mathrm{SP4} \\
\mathrm{adj0}_{t-4}^{\mathrm{L}}, & \text{if } t > \mathrm{SP4}
\end{cases}, \qquad
\mathrm{adj4}_{t}^{\mathrm{H}} =
\begin{cases}
0, & \text{if } t \leq \mathrm{SP4} \\
\mathrm{adj0}_{t-4}^{\mathrm{H}}, & \text{if } t > \mathrm{SP4}
\end{cases}, && \forall\, t \in \mathcal{T}_\mathrm{SP}, \label{eq:neso_soe_mana:adj0_adj4}\\
& \mathrm{left}_{t}^{\mathrm{L}} = 0,\qquad
  \mathrm{left}_{t}^{\mathrm{H}} = 0, && \forall t \in \mathcal{T}_{\mathrm{EFA}}^{\mathrm{SP1}}, \label{eq:neso_soe_mana:left_sp0}\\
& \mathrm{left}_{t}^{\mathrm{L}} = \mathrm{left}_{t-1}^{\mathrm{L}} + \mathrm{FRE}_{t-1}^{\mathrm{L}} - \mathrm{RER}_{t-1}^{\mathrm{L}},
&& \forall t\in\mathcal{T}_{\mathrm{EFA}}^{\mathrm{SP2-8}}, \label{eq:neso_soe_mana:left_cal_L}\\
& \mathrm{left}_{t}^{\mathrm{H}} = \mathrm{left}_{t-1}^{\mathrm{H}} + \mathrm{FRE}_{t-1}^{\mathrm{H}} - \mathrm{RER}_{t-1}^{\mathrm{H}},
&& \forall t\in\mathcal{T}_{\mathrm{EFA}}^{\mathrm{SP2-8}}.\label{eq:neso_soe_mana:left_cal_H}
\end{align}
\end{subequations}

Binary decision variables $b_t^{\mathrm{H}}$ and $b_t^{\mathrm{L}}$ are defined by constraints~\eqref{eq:neso_soe_mana:bL_def}–\eqref{eq:neso_soe_mana:bH_def} to indicate the provision of high-frequency and low-frequency response in timestep $t$, respectively. The possible combinations \( (b_t^{\mathrm{H}}, b_t^{\mathrm{L}}) = (0,0) \), \( (1,0) \) or \( (0,1) \), and \( (1,1) \) correspond to no contract, single-directional contracts, and a bi-directional contract, respectively. Based on these configurations, constraints~\eqref{eq:neso_soe_mana:reserve_L}--\eqref{eq:neso_soe_mana:reserve_H_bound2} jointly determine the reserve amounts \(r_t^{\mathrm{L}}\) and \(r_t^{\mathrm{H}}\) based on the contract type.  
When \( (b_t^{\mathrm{H}}, b_t^{\mathrm{L}}) = (1,1) \) (i.e., bi-directional contract), the right-hand sides of \eqref{eq:neso_soe_mana:reserve_L}--\eqref{eq:neso_soe_mana:reserve_H} compute the 20\% energy recovery reserve requirement that must be allocated in each direction in MWh. Following the discussion in~\ref{app:neso_soe}, the reserve factors are set to \(10\%\), \(20\%\), and \(40\%\) of the offered service power \(p_t^{\mathrm{DC}}\), \(p_t^{\mathrm{DM}}\), and \(p_t^{\mathrm{DR}}\) in MW, respectively. When the contract is single-directional, constraints \eqref{eq:neso_soe_mana:reserve_L_bound1}--\eqref{eq:neso_soe_mana:reserve_H_bound2} force the entire reserve to the non-contracted direction, while the contracted-direction reserve is fixed at zero. If no service is contracted (\(b_t^{\mathrm{H}}=b_t^{\mathrm{L}}=0\)), both reserve variables are set to zero. The product \(b_t^{\mathrm{H}}b_t^{\mathrm{L}}\) is linearised using an auxiliary binary variable that equals one only when both \(b_t^{\mathrm{H}}\) and \(b_t^{\mathrm{L}}\) are one.


Constraints~\eqref{eq:neso_soe_mana:L_pmax} and \eqref{eq:neso_soe_mana:H_pmax} ensure that the aggregate charging and discharging power, including baseline operation, contracted services, and reserves, remains within the rated power capacity. Next, constraints~\eqref{eq:neso_soe_mana:CREV_L} and \eqref{eq:neso_soe_mana:CREV_H} compute the contracted response energy volume for low-frequency and high-frequency services in each 4-hour electricity forward agreement block. In each case, the contracted response energy volume is calculated by converting the contracted power of the DC, DM, and DR products into energy through multiplication by their respective delivery durations (see Table~\ref{tab:dfr_technical_requirements}). The energy recovery requirement is subsequently set to 20\% of the contracted response energy volume, expressed in MWh, as specified in~\eqref{eq:neso_soe_mana:ER}.

The frequency response energy \(\mathrm{FRE}_t^{\mathrm{L}}\) and \(\mathrm{FRE}_t^{\mathrm{H}}\), defined in \eqref{eq:neso_soe_mana:FRE_L_drf} and \eqref{eq:neso_soe_mana:FRE_H_drf}, denote the energy delivered by low- and high-frequency response, respectively, during settlement period $t$. These quantities represent the energy delivered during one settlement period, calculated by applying the frequency-deviation activation factors $\alpha$ to the contracted response capacity over the settlement period. The activation functions that map frequency deviations to $\alpha$ are shown in Fig.~\ref{fig:fre_act}. Constraints~\eqref{eq:neso_soe_mana:REV_L_bound} and \eqref{eq:neso_soe_mana:REV_H_bound} impose upper bounds on \(\mathrm{FRE}_t^{\mathrm{L}}\) and \(\mathrm{FRE}_t^{\mathrm{H}}\), ensuring that the energy delivered within each settlement period remains within the available 30-minute energy capacity of the asset.

Although the optimisation is formulated under perfect foresight, uncertainty in actual service activation is represented through conservative operational constraints across different contracted service periods. To this end, constraints~\eqref{eq:neso_soe_mana:offset_H} and \eqref{eq:neso_soe_mana:offset_L} require baseline charging and discharging actions to account for the frequency response energy delivered within each contracted service period. Specifically, a minimum level of baseline discharging (charging) power is scheduled in proportion to the realised high-frequency (low-frequency) response energy within each electricity forward agreement block, scaled by a factor $\delta$. This ensures that a portion of the energy exchanged through service activation is compensated for within the same electricity forward agreement block, enabling continuous service delivery and compliance with minimum SoE requirements at the start of each settlement period for all contracted service periods.

Constraint~\eqref{eq:neso_soe_mana:MSER_START} then initialises the minimum SoE requirement at the start of each electricity forward agreement block as equal to the contracted response energy volume, reflecting the energy needed to fulfil contracted obligations. Note that in NESO's terminology \cite{NESO_SOE_monitor}, the minimum state-of-energy requirement (MSER) is expressed in MWh, rather than as a normalised SoE value between 0 and 1. To allow flexibility in SoE management, constraint~\eqref{eq:neso_soe_mana:MSER_allSP} permits the minimum SoE requirement to be updated in later settlement periods within the same block, while ensuring this never exceeds the contracted response energy volume. Constraints~\eqref{eq:neso_soe_mana:MSER_cal_L} and \eqref{eq:neso_soe_mana:MSER_cal_H} update the minimum SoE requirement across consecutive settlement periods, with its evolution governed by the delivered energy and the implemented energy recovery adjustment volume~\cite{NESO_Service_Terms}. Constraints~\eqref{eq:neso_soe_mana:soc_L_bound} and \eqref{eq:neso_soe_mana:soc_H_bound} bound the battery SoE at the start of each settlement period, ensuring sufficient upward and downward SoE margins throughout the contracted service period.

Next, we explain how to determine the energy adjustment volume $\mathrm{adj4}_t$ resulting from DFR service delivery. The energy margin variables \(\mathrm{MG}_t^{\mathrm{L}}\) and \(\mathrm{MG}_t^{\mathrm{H}}\), defined in \eqref{eq:neso_soe_mana:mg_cal_L} and \eqref{eq:neso_soe_mana:mg_cal_H}, measure the deviation of the current SoE from the contracted response energy volume bounds for low- and high-frequency services, respectively. Rather than enforcing a full return to the initial SoE level, the mechanism uses the energy margin variables to compute the minimum corrective energy adjustment required to restore the contracted response energy volume feasibility. Equations \eqref{eq:neso_soe_mana:RER_cal_L} and \eqref{eq:neso_soe_mana:RER_cal_H} define the required energy recovery $\mathrm{RER}_t$ as the recovery obligation arising from service delivery \(\mathrm{FRE}_t\) and any carry-over requirement (\(\mathrm{left}_t\)), capped by the corresponding energy recovery limits per settlement period (see NESO guidance from \cite{NESO_SOE_monitor}). Constraints~\eqref{eq:neso_soe_mana:ABSR_cal_L} and~\eqref{eq:neso_soe_mana:ABSR_cal_H} define the absorbed energy $\mathrm{ABS}_t$ as the portion of the required energy recovery that can be absorbed by the available SoE margin, thereby limiting the amount of recovery that must be scheduled through baseline adjustment (see Fig.~\ref{fig:neso_case2}) \cite{NESO_new_explain}. Equation~\eqref{eq:neso_soe_mana:ad0_def} defines $\mathrm{adj0}_t$ as the remaining recovery requirement after accounting for the absorbed energy. This design ensures that SoE corrections are limited to what is necessary for contract compliance, thereby avoiding unnecessary energy recovery and preserving operational flexibility. Constraint \eqref{eq:neso_soe_mana:adj0_adj4} enforces a 4-settlement-period implementation delay (see Section~\ref{sec:bg:soe}), such that the energy adjustment scheduled in timestep $t$, $\mathrm{adj0}_t$, only becomes effective at $t+4$ as $\mathrm{adj4}_{t+4}$.


Finally, the carry-over term `$\mathrm{left}_t$' is introduced to track unrecovered energy. At the start of each electricity forward agreement block, $\mathrm{left}_t$ is reset to zero by \eqref{eq:neso_soe_mana:left_sp0}. During the block, \eqref{eq:neso_soe_mana:left_cal_L} and \eqref{eq:neso_soe_mana:left_cal_H} update $\mathrm{left}_t$ in each settlement period. The carry-over term is calculated as the difference between the frequency response energy delivered in the previous settlement period and the corresponding required energy recovery. It represents any recovery shortfall that cannot be met within a single settlement period due to energy recovery limits.

\subsection{High fidelity battery model}\label{sec:mpc:hfbm}
After solving the first-stage optimisation model over the optimisation horizon, the resulting AC power dispatch $p_t^{\mathrm{tot}}$ with time step size $\Delta t$ is obtained. This dispatch is subsequently simulated using SimSES~\cite{SimSES_paper} to evaluate battery operation and degradation. SimSES is an open-source simulation environment for the techno-economic analysis of battery energy storage systems, supporting a wide range of system designs and control strategies, including grid-scale BESS applications. The framework has been widely applied in the literature to analyse battery operation and degradation in both mobility and stationary battery systems~\cite{MPC_battery,8916568,TEPE2022118351,PARLIKAR2023120541}, and the complete code base is publicly available \cite{SimSES_paper}. The parametrisation adopted for the BESS model in this study is summarised in Table~\ref{tab:simses_parameter}.

\begin{table}[t]
\centering
\small
\caption{Model parameters and values for the BESS simulation in SimSES \cite{SimSES_paper}.}
\begin{tabular}{p{5cm} p{8cm}}
\toprule
\midrule
\textbf{Parameter} & \textbf{Value} \\
\midrule
Cell type & Graphite LiFePO$_4$ (Sony US26650FTC1) \cite{cyc_ageing_LFP,SimSES_paper,calen_ageing_LFP} \\
Rated energy capacity (MWh) & 5 \\
Rated power (MW) & 5 \\
Battery model & Equivalent circuit model \cite{SimSES_paper} \\
Battery degradation model & Semi-empirical, calendar and cycle ageing \cite{cyc_ageing_LFP,calen_ageing_LFP}\\
AC/DC Converter & Notton et al.\ \cite{NOTTON_converter} \\
\midrule
\bottomrule
\end{tabular}
\label{tab:simses_parameter}
\end{table}

\section{Case study}\label{sec:case_study}
In this section, we present simulation results obtained using the proposed model predictive control framework, within which six ageing-aware optimisation models are implemented to evaluate battery lifetime profitability and operational strategies in the GB electricity and frequency response markets. Market data for the year 2024 are used as the reference. System frequency data for the GB grid are publicly available~\cite{neso_system_frequency}, while DFR clearing prices are sourced from the NESO~\cite{neso_eac_auction_results}. Intra-day electricity prices, represented by the volume-weighted average price, are obtained from EPEX SPOT~\cite{epexspot_gbdata}. All additional assumptions and model parameters are detailed in~\ref{sec:para_setting} and Table~\ref{tab:model_params}. The simulation horizon is set to 20 years, with termination occurring earlier if the battery reaches the EoL criterion of 80\% SoH. All simulations are executed using Gurobi Optimizer version 13.0.2 on a high-performance computing server equipped with dual AMD EPYC 7453 processors (112 threads) and 256 GB of RAM. The following sections present the main analyses. Section~\ref{sec:cs:ageing_profit} investigates the impact of different ageing models on BESS lifetime profitability. Next, Section~\ref{sec:cs:dr_profit} examines the effect of discount rates. Finally, Section~\ref{sec:cs:ageing_service} analyses how ageing models and past degradation influence the optimal service stacking strategy.

\subsection{Impact of different ageing models on economic outcomes}\label{sec:cs:ageing_profit}
\begin{figure}[t]
  \centering
  \includegraphics[width=0.6\textwidth]{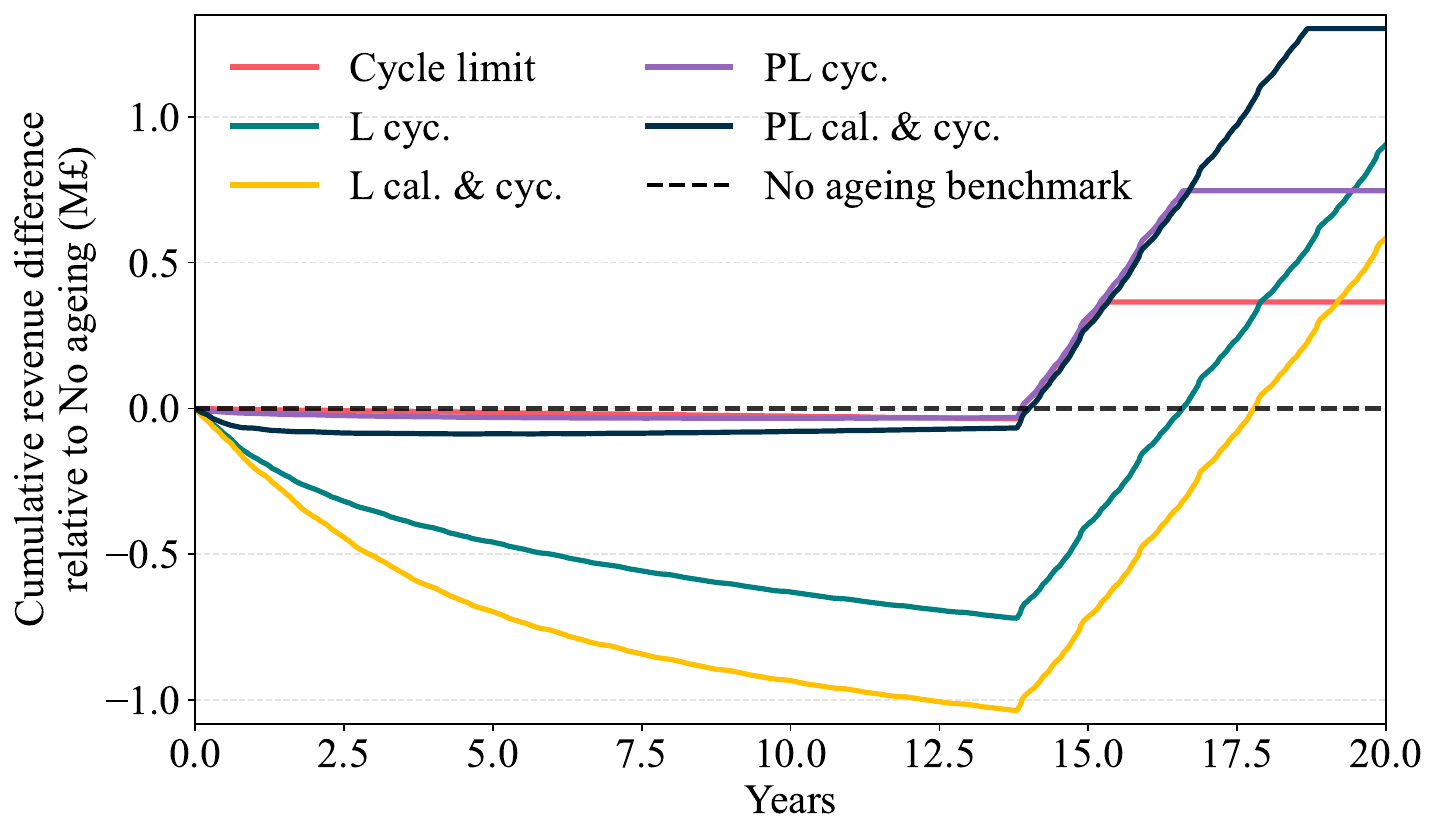}
  \caption{Cumulative revenue difference of ageing-aware optimisation results relative to the \textit{No ageing} benchmark over the project lifetime (20 years). The dashed horizontal line indicates zero difference, corresponding to equal performance with the \textit{No ageing} model.}
  \label{fig:revenue_difference}
\end{figure}

\begin{figure}[t]
  \centering
  \includegraphics[width=\textwidth]{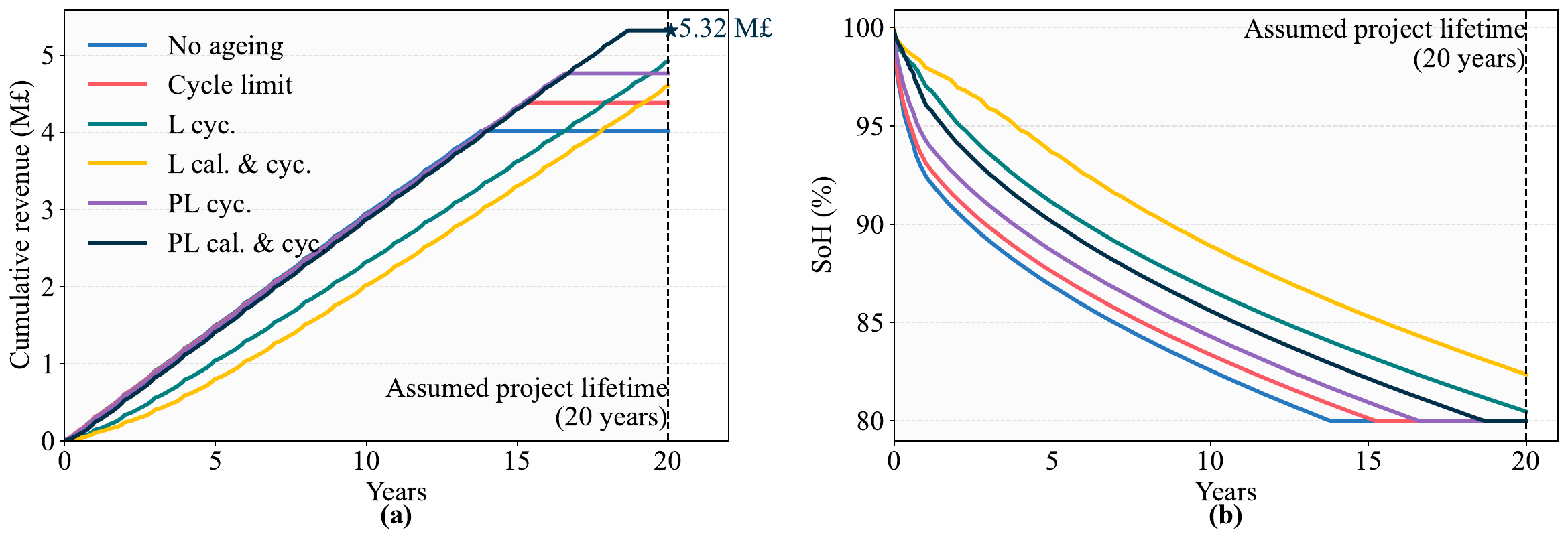}
  \caption{Cumulative lifetime economic performance and degradation trajectories under different ageing-aware models: (a) cumulative revenue; (b) state-of-health (SoH) evolution over a 20-year horizon.}
  \label{fig:lifetime_revenue}
\end{figure}

\begin{figure}[t]
  \centering
  \includegraphics[width=\textwidth]{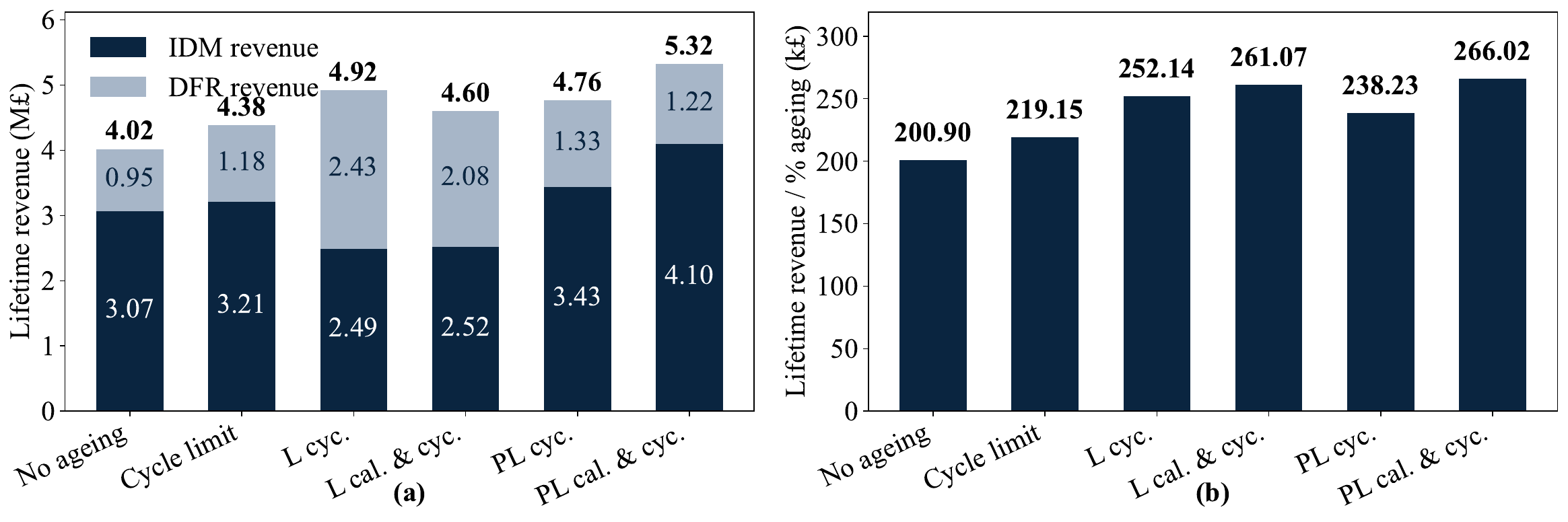}
  \caption{Lifetime revenue composition and economic return per degradation under different ageing-aware models: (a) total revenue breakdown by market; (b) lifetime revenue per \% ageing.}
  \label{fig:rev_detail_per_plot}
\end{figure}

The case study investigates simulated BESS operation across six battery ageing model formulations, resulting in different operating strategies with differing lifetime economic performance and degradation characteristics. The results are presented in Table~\ref{tab:annual_pro_results} and Figs.~\ref{fig:revenue_difference}--\ref{fig:rev_detail_per_plot}.

\begin{table}[t]
\centering
\small
\caption{Lifetime cumulative simulation results for different ageing-aware models.}
\label{tab:annual_pro_results}
\setlength{\tabcolsep}{6pt}
\begin{threeparttable}
\begin{tabular}{l c c c c c c}
\toprule
\midrule
 \textbf{Ageing model}
& \textbf{\shortstack{$\Pi^{\mathrm{E}}$ \\ (M\pounds)}}
& \textbf{\shortstack{$\Pi^{\mathrm{DFR}}$ \\ (M\pounds)}}
& \textbf{\shortstack{Revenue \\ (M\pounds)}}
& \textbf{\shortstack{EoS SoH \\ (\%)}}
& \textbf{\shortstack{Time to EoL \\ 80\% SoH (years)}}
& \textbf{\shortstack{Revenue per \% \\ of ageing   (k\pounds)}}\\
\midrule
\textit{No ageing} & 3.07 & 0.95 & 4.02 & 80.00 & 13.80 & 200.90 \\
\textit{Cycle limit} & 3.21 & 1.18 & 4.38 & 80.00 & 15.23 & 219.15 \\
\textit{L cyc.} & 2.49 & 2.43 & 4.92 & 80.48 & N/A & 252.14 \\
\textit{L cal. \& cyc.} & 2.52 & 2.08 & 4.60 & 82.37 & N/A & 261.07 \\
\textit{PL cyc.} & 3.43 & 1.33 & 4.76 & 80.00 & 16.59 & 238.23 \\
\textit{PL cal. \& cyc.} & 4.10 & 1.22 & 5.32 & 80.00 & 18.67 & 266.02 \\
  \midrule
\bottomrule
\end{tabular}
\begin{tablenotes}[flushleft]
\footnotesize
\item \textit{Notes:} 
Revenue is defined as the sum of electricity market revenue $\Pi^{\mathrm{E}}$ and DFR market revenue $\Pi^{\mathrm{DFR}}$. The EoS SoH denotes the battery SoH at the end-of-simulation (EoS) horizon, obtained from SimSES. Revenue per \% of ageing is defined as 
$(\Pi^{\mathrm{E}}+\Pi^{\mathrm{DFR}})/(100-\text{EoS SoH})$.
\end{tablenotes}
\end{threeparttable}
\end{table}

Fig.~\ref{fig:revenue_difference} illustrates the evolution of cumulative revenue differences for five ageing-aware models relative to the \textit{No ageing} benchmark over the project lifetime. The \textit{No ageing} strategy prioritises short-term revenue by allowing unconstrained battery operation, thereby fully exploiting market opportunities without accounting for degradation. This leads to higher revenues in the early years, reflected by negative revenue differences across all ageing-aware strategies. However, this advantage diminishes as accelerated degradation drives the battery to its EoL threshold at approximately 13.8 years under the \textit{No ageing} strategy. Beyond this point, ageing-aware strategies continue to generate revenue due to their extended operational lifetimes, resulting in a reversal of revenue differences. As a consequence, despite its initial performance advantage, the \textit{No ageing} model achieves a total lifetime revenue of £4.02~M, which is approximately 24\% lower than that of the best-performing ageing-aware strategy, \textit{PL cal. \& cyc.} (£5.32~M). This demonstrates the significant economic value of incorporating degradation-aware decision-making.

Fig.~\ref{fig:lifetime_revenue} and Table~\ref{tab:annual_pro_results} show that different operational strategies lead to variations in operational lifetime. The \textit{No ageing} model reaches the EoL threshold prematurely, whereas ageing-aware strategies explicitly account for degradation and extend the operational lifetime to varying degrees. For example, the \textit{Cycle limit} and \textit{PL cyc.} models extend the time to EoL to 15.2 and 16.6 years, respectively, while the \textit{PL cal. \& cyc.} model further extends the lifetime to 18.7 years. 

However, lifetime extension alone does not guarantee improved economic performance. Aggressive strategies (i.e., \textit{No ageing} and \textit{Cycle limit}) prioritise short-term gains but fail to preserve sufficient lifetime, resulting in limited lifetime profitability. This is also reflected in Fig.~\ref{fig:lifetime_revenue}(b), where these strategies exhibit faster SoH degradation and reach the EoL threshold significantly earlier. In contrast, overly conservative strategies (i.e., \textit{L cyc.} and \textit{L cal. \& cyc.}) achieve extended operational lifetimes but significantly restrict market participation. These linear models maintain SoH above the EoL threshold over the 20-year horizon, as shown in Fig.~\ref{fig:lifetime_revenue}(b), indicating effective degradation control under strong penalisation. However, this reflects an overly conservative operating regime that limits asset utilisation and suppresses revenue accumulation. As a result, despite sustaining operation over the full project lifetime, their cumulative revenues remain constrained. Among all strategies, the \textit{PL cal. \& cyc.} model achieves the highest lifetime revenue (£5.32~M), representing a 32.34\% increase over the \textit{No ageing} model and a 15.65\% improvement compared to the best non-piecewise model (\textit{L cal. \& cyc.}). This highlights the importance of high-fidelity degradation modelling in capturing nonlinear dependencies on operating conditions. Overall, this shows that lifetime revenue is maximised through a balance between degradation and value extraction over time.

To further examine how economic value is extracted from battery degradation, Fig.~\ref{fig:rev_detail_per_plot}(b) presents the lifetime revenue per \% of capacity loss, capturing the value generated from each unit of degradation. The \textit{PL cal. \& cyc.} model achieves both the highest total revenue (£5.32~M) and the highest economic value per unit of degradation (£266.02~k), indicating efficient utilisation of the asset. By contrast, the \textit{L cal. \& cyc.} model exhibits a relatively high revenue per \% of ageing despite lower total revenues. This is driven by a shift in optimal service stacking under strong degradation penalisation, leading to operating strategies that reduce degradation while prioritising higher-value utilisation of each unit of degradation. As a result, more economic value is extracted per unit of ageing, even though total revenue remains constrained. Conversely, the \textit{PL cyc.} model adopts more aggressive utilisation, achieving higher total revenue at the expense of faster degradation, leading to a lower revenue per unit of ageing.

The impact of calendar ageing is strongly dependent on ageing-aware model formulation. Under the simple linear formulation, incorporating calendar ageing increases the overall degradation cost, reinforcing an already conservative operating regime and reducing lifetime revenue over 20 years by up to 6.5\%. In this case, degradation is effectively over-penalised, and the addition of calendar ageing further suppresses asset utilisation, leading to lower revenue over the 20-year horizon. In contrast, under the sophisticated piecewise linear formulation, incorporating calendar ageing leads to an 11.8\% increase in lifetime revenue (from \textit{PL cyc.} to \textit{PL cal. \& cyc.}). The piecewise linear formulation avoids the systematic over-penalisation of degradation present in the linear model, allowing the inclusion of calendar ageing to influence operational decisions without inducing excessive conservativeness. This results in a longer effective operational lifetime (from 16.59 to 18.67 years) and improved utilisation of the battery.

However, Fig.~\ref{fig:rev_detail_per_plot}(b) shows that incorporating calendar ageing \textit{consistently} increases the revenue generated per unit of battery degradation across both modelling approaches, by approximately 3.5\% under the linear formulation and 11.7\% under the piecewise linear formulation. This is because the inclusion of calendar ageing alters the optimal operating regime, as degradation becomes dependent not only on energy throughput but also on SoC. As a result, the model adjusts its operation towards conditions that mitigate degradation, such as maintaining lower SoC levels and shifting the service mix, while still capturing revenue opportunities. Consequently, each unit of degradation is associated with higher revenue generation, leading to improved economic returns per unit of capacity loss. Overall, these results indicate that the impact of calendar ageing on total revenue is conditional on model fidelity, but it consistently improves the value extracted per unit of degradation.

These differences in operational strategies are ultimately reflected in market participation patterns. The detailed impact of degradation modelling on optimal service stacking is discussed later in Section~\ref{sec:cs:ageing_service}. As shown in Fig.~\ref{fig:rev_detail_per_plot}(a), revenue from IDM dominates that from DFR, indicating that IDM is the primary source of economic value. However, ageing-aware strategies lead to different service stacking decisions, which determine how the battery participates across markets and ultimately drive total revenue. More conservative models rely more heavily on DFR services, particularly short-duration products, which involve lower energy throughput and reduced degradation rates, but also limit access to higher-value IDM opportunities. In contrast, less conservative strategies allocate more operation to IDM, which increases revenue potential but also results in higher degradation. This demonstrates that market participation is a key mechanism through which degradation effects are translated into economic value.

\subsection{Impact of different discount rates on economic outcomes}\label{sec:cs:dr_profit}
In this section, we investigate how discount rates influence the economic outcomes of battery operation under different ageing-aware strategies. By incorporating discounted lifetime revenue, we examine the trade-off between short-term revenue generation and long-term asset utilisation, while accounting for degradation effects. The lifetime revenue is calculated by aggregating daily revenues from energy arbitrage $\Pi^{\mathrm{E}}$ and dynamic frequency response services $\Pi^{\mathrm{DFR}}$, and discounting future cash flows using an annual discount rate. Revenue accumulation is terminated once the battery reaches the EoL threshold, defined as a SoH below 80\%.

\begin{figure}[t]
  \centering
  \includegraphics[width=\textwidth]{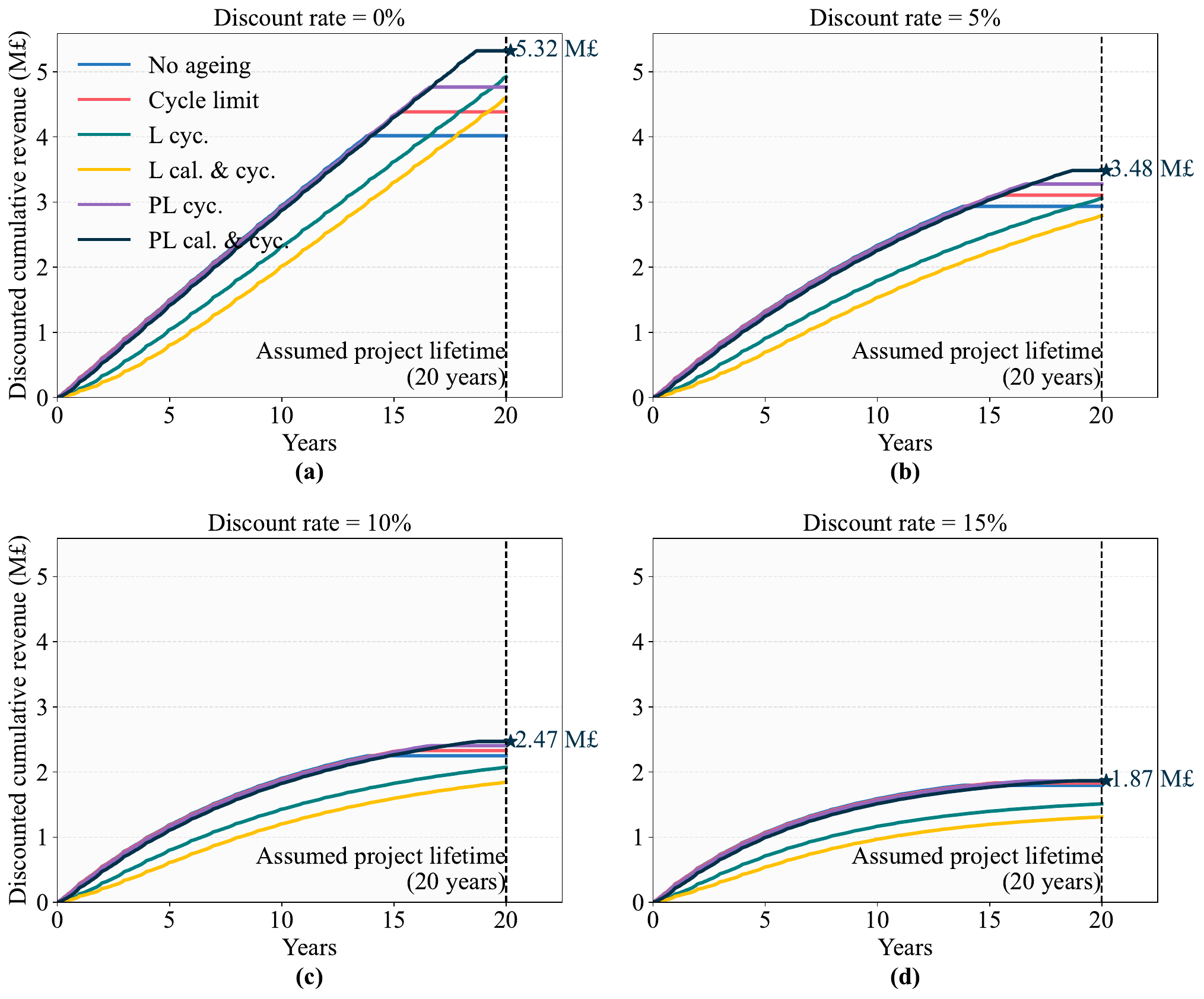}
  \caption{Evolution of cumulative discounted revenue assuming an annual discount rate of 0, 5\%, 10\% and 15\%. Black line represents assumed project lifetime of 20 years.}
  \label{fig:discount_1}
\end{figure}

\begin{figure}[t]
  \centering
  \includegraphics[width=0.7\textwidth]{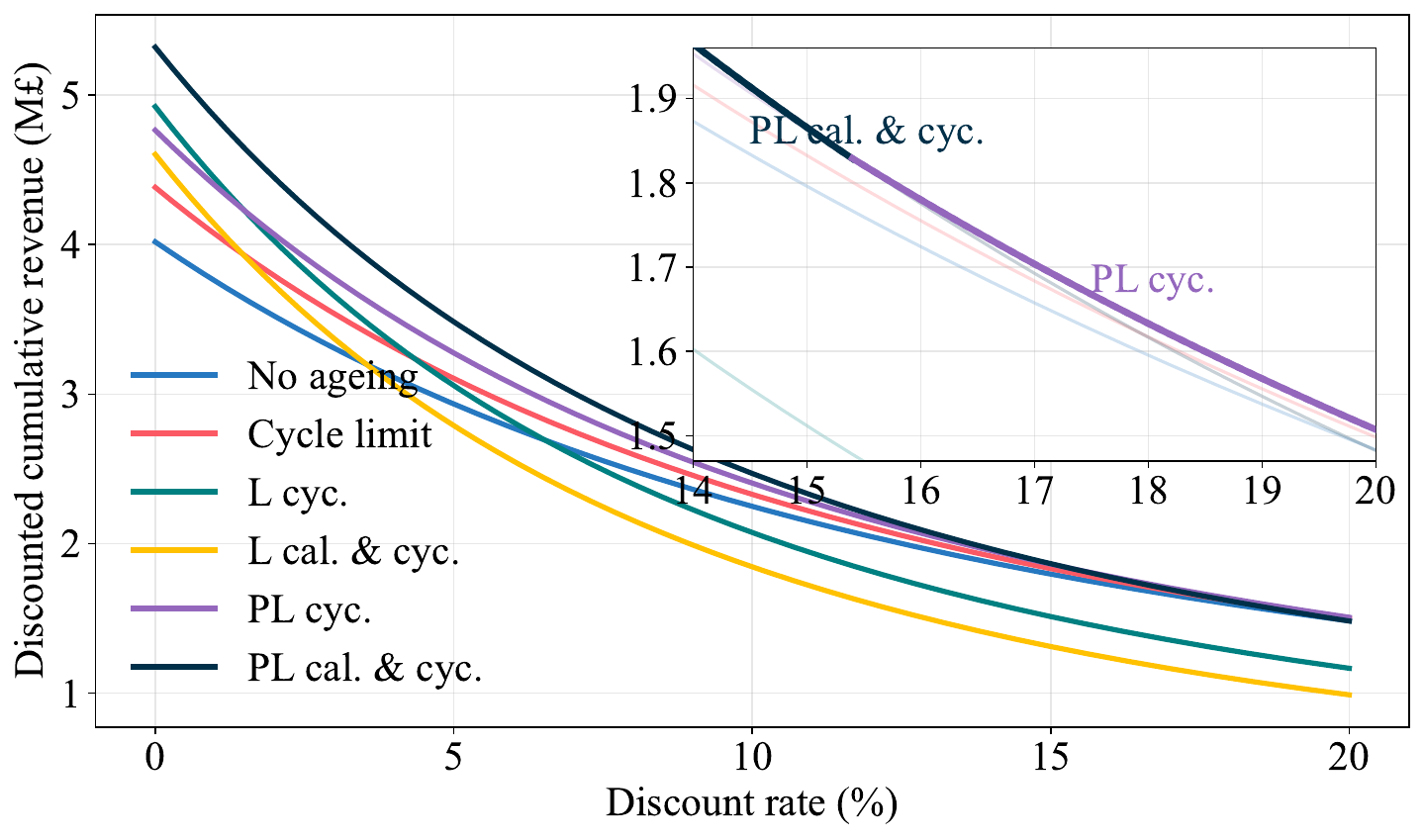}
  \caption{Impact of discount rate on discounted lifetime revenue under different ageing-aware models. Results are evaluated over discount rates from 0\% to 20\% with a step size of 0.1\%. The inset enlarges the high-discount-rate region, where the strategy achieving the maximum discounted lifetime revenue changes with the discount rate.}
  \label{fig:discount_2}
\end{figure}

As shown in Fig.~\ref{fig:discount_1}, for discount rates of 0\%, 5\%, 10\%, and 15\%, the \textit{PL cal. \& cyc.} model consistently remains the optimal strategy, achieving £5.32 M, £3.48 M, £2.47 M, and £1.87 M, respectively. However, its advantage over the second-best strategy gradually diminishes as the discount rate increases. This trend reflects the reduced economic value of long-term revenue under higher discount rates. The \textit{PL cal. \& cyc.} model derives its advantage from both extending battery lifetime and maintaining effective utilisation, enabling sustained revenue generation without overly conservative operation. As the discount rate increases, these long-term benefits are progressively discounted, weakening its relative advantage. Consequently, strategies that generate revenue earlier become increasingly competitive under higher discount rates.

To further investigate this effect, a sensitivity analysis was conducted over discount rates ranging from 0\% to 20\% with a step size of 0.1\%, as shown in Fig.~\ref{fig:discount_2}. The results reveal clear strategy transitions as the discount rate increases. Specifically, the \textit{PL cal. \& cyc.} model is outperformed by the \textit{PL cyc.} model at 15.4\%, by the \textit{Cycle limit} model at 18\%, and ultimately by the \textit{No ageing} strategy at 19.9\%. These transitions indicate that the optimal operational strategy is sensitive to the valuation of future revenues. At low discount rates, future revenues are valued more strongly, favouring strategies that preserve battery capacity and defer degradation. As the discount rate increases, future revenues are progressively discounted, reducing the economic penalty associated with accelerated degradation. Consequently, strategies that generate higher near-term revenues become increasingly favourable, even if they lead to faster capacity loss.

Fig.~\ref{fig:discount_2} further shows that \textit{PL cal. \& cyc.} remains the optimal strategy over a wide range of discount rates, indicating that accounting for both calendar and cycling ageing is economically favourable under moderate discounting conditions. This result is expected because real batteries lose capacity not only through cycling, but also through SoC- and time-dependent calendar ageing. A cycle-only model ignores this loss and can therefore overvalue the benefit of saving battery capacity for future use. At higher discount rates, however, the optimal strategy transitions to \textit{PL cyc.}, which allows more aggressive battery utilisation by imposing a lower effective degradation penalty. This becomes favourable when future revenues and future degradation impacts are strongly discounted. The small performance difference between \textit{PL cal. \& cyc.} and \textit{PL cyc.} at high discount rates suggests that the theoretical benefit of including calendar ageing is weakened in the numerical results, likely because high discounting reduces the economic value of long-term degradation effects and the piecewise-linear approximation may not fully preserve the nonlinear behaviour of calendar ageing.

\subsection{Impact of ageing-aware modelling and degradation state on service allocation}\label{sec:cs:ageing_service}
\begin{figure}[htbp] 
  \centering
  \includegraphics[width=\columnwidth]{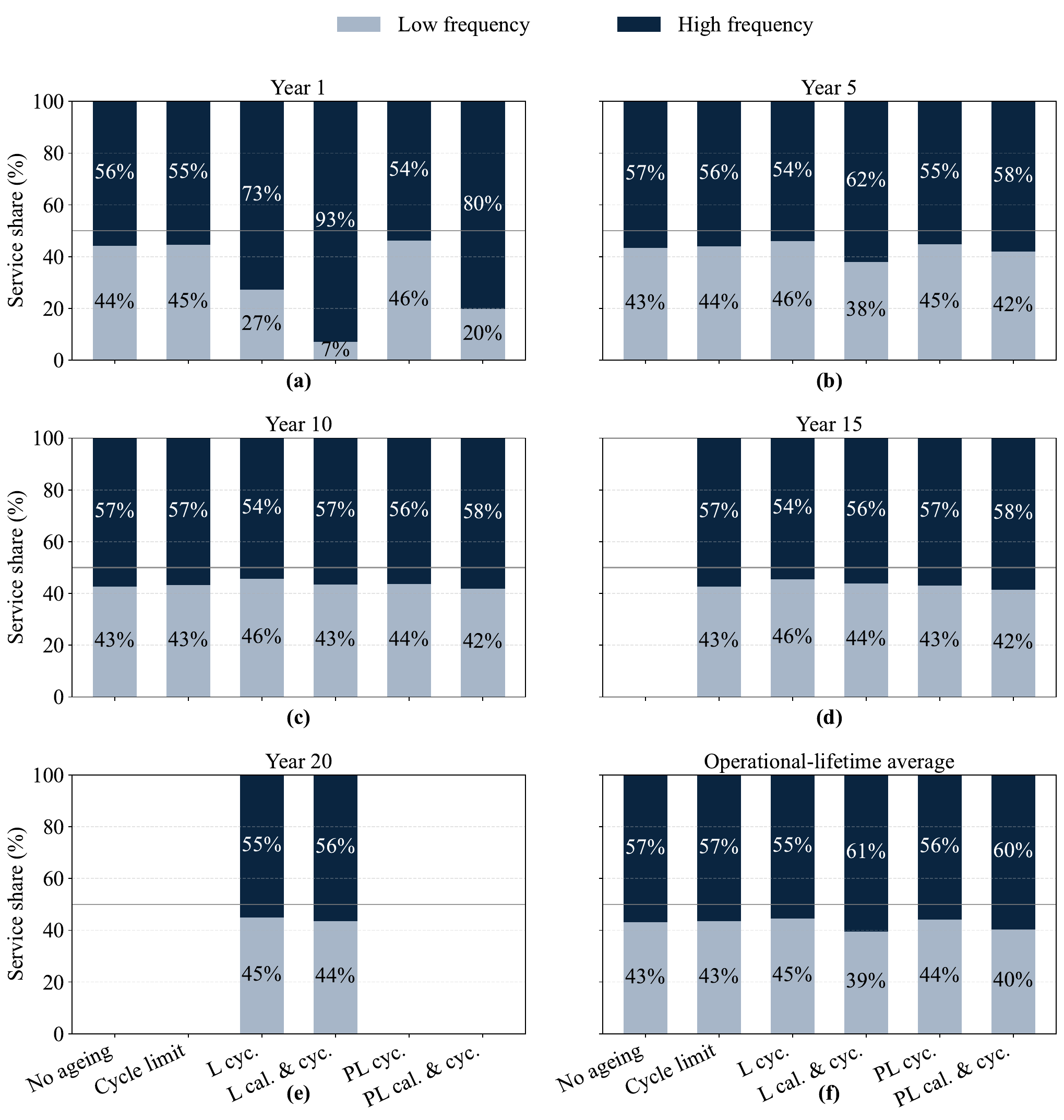}
  \caption{Service allocation and temporal evolution across high- and low-frequency DFR under different ageing-aware models: (a–e) representative years; (f) lifetime average. The absence of bars indicates that the battery has reached end-of-life (EoL) in the corresponding year.}\label{fig:service_H_L}
\end{figure}

\begin{figure}[htbp] 
  \centering
  \includegraphics[width=\columnwidth]{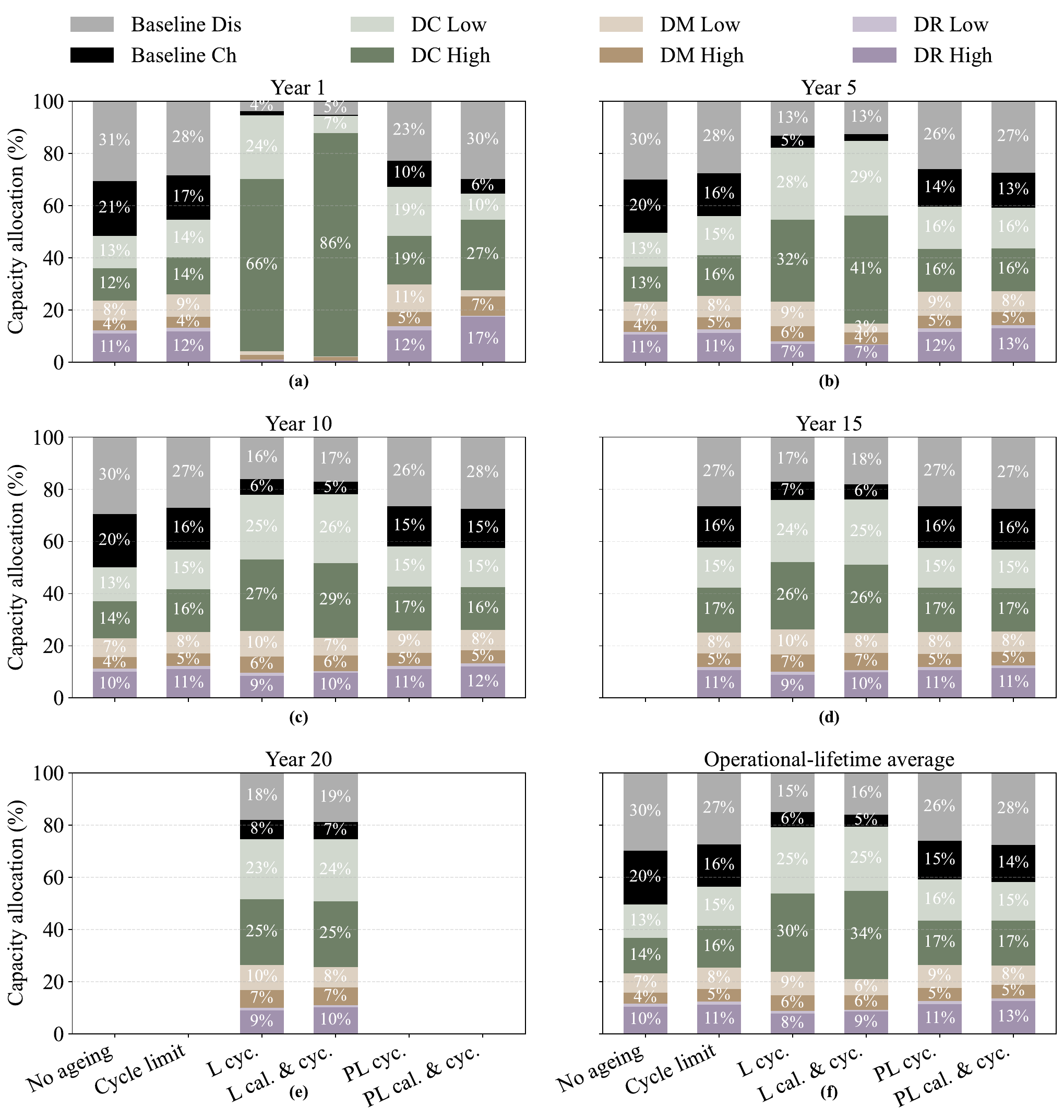}
  \caption{Capacity allocation across services under different ageing-aware models: (a–e) representative years; (f) lifetime average. The absence of bars indicates that the battery has reached end-of-life (EoL) in the corresponding year.}\label{fig:service_year}
\end{figure}

\begin{figure}[htbp] 
  \centering
  \includegraphics[width=\columnwidth]{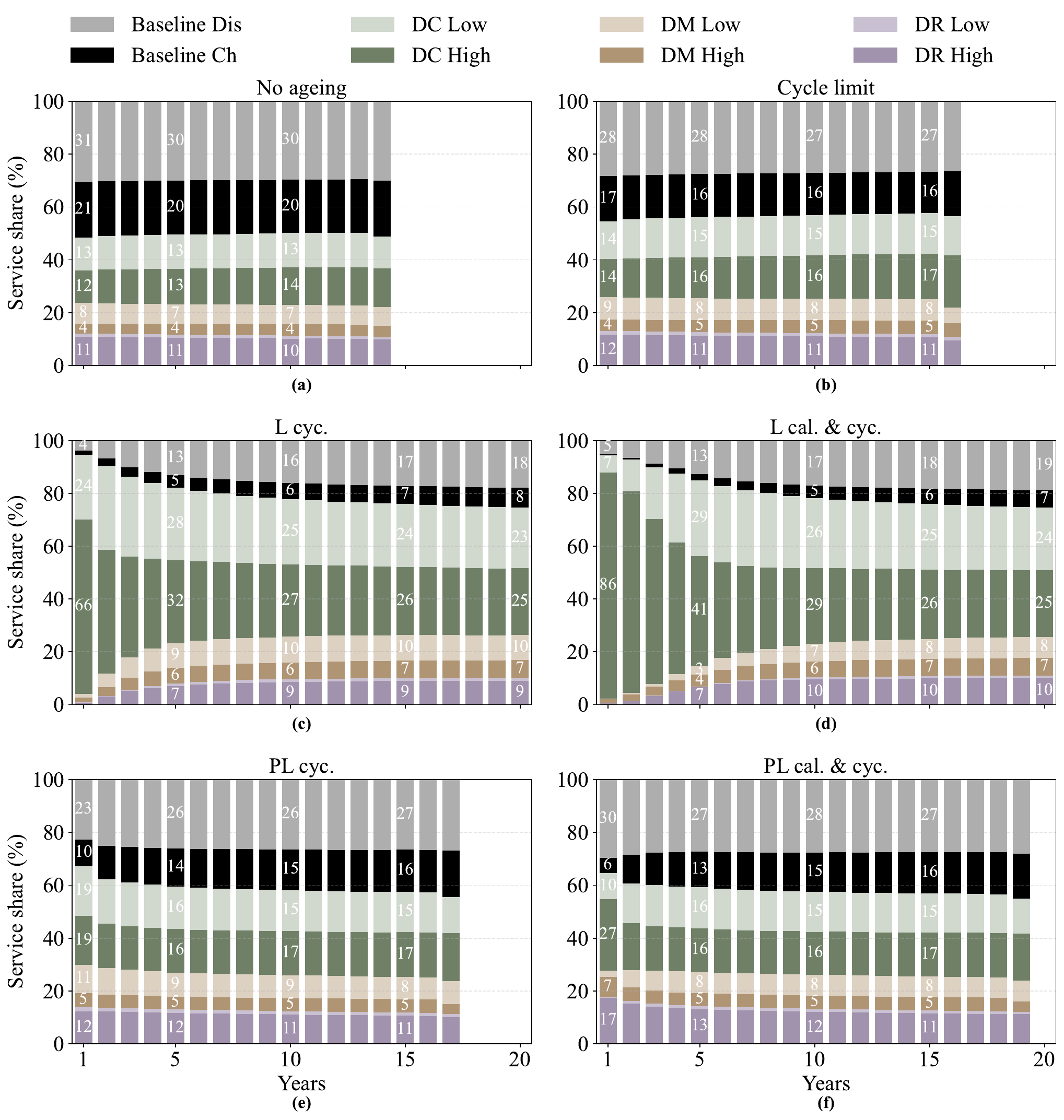}
  \caption{Temporal evolution of capacity allocation across services under different ageing-aware models. The absence of bars indicates that the battery has reached end-of-life (EoL) in the corresponding year.}\label{fig:service_ageing}
\end{figure}

To complement the economic analysis, we next examine how different ageing models influence optimal service stacking over the lifetime simulation horizon. Figs.~\ref{fig:service_H_L}--\ref{fig:service_ageing} illustrate the annual average service distribution across all ageing-aware formulations.

Fig.~\ref{fig:service_H_L} shows a clear and consistent preference for high-frequency services over low-frequency services across all models. This pattern is mainly driven by the interaction between DFR provision and baseline operation. Specifically, high-frequency service activation charges the battery, which can subsequently support baseline discharge and thereby create additional revenue opportunities from energy arbitrage. By contrast, low-frequency services require the battery to deliver energy to the grid, reducing its SoE and increasing the need for subsequent energy recovery through baseline charging, which incurs additional cost. This tendency is particularly pronounced in calendar ageing-aware strategies, as high-frequency service provision not only enhances revenue by facilitating subsequent baseline discharge, but also helps maintain lower average SoC levels, thereby mitigating calendar degradation. Consequently, DFR high-frequency services are consistently preferred over low-frequency across all strategies.

Beyond the clear preference for high-frequency over low-frequency products, Figs.~\ref{fig:service_year} and \ref{fig:service_ageing} reveal strategy-dependent service-stacking patterns across the full set of products. The \textit{No ageing} and \textit{Cycle limit} strategies exhibit relatively diversified service portfolios, with operation distributed fairly evenly across services. Nevertheless, baseline discharging remains the largest single component, accounting for average lifetime shares of approximately 30\% and 28\%, respectively. In contrast, the L-based models (i.e., \textit{L cyc.} and \textit{L cal. \& cyc.}) exhibit highly concentrated service allocation profiles, with a dominant share assigned to DCH services and, to a lesser extent, DCL services and baseline discharge. These DC services are characterised by short duration and low energy throughput, and are therefore favoured by strategies with stronger degradation penalisation. This concentration reflects the conservative behaviour induced by the linear degradation representation, which assigns a relatively high degradation cost to energy-intensive operation and therefore shifts the optimisation toward low-throughput services. The piecewise linear formulations provide a more balanced allocation. Both \textit{PL cyc.} and \textit{PL cal. \& cyc.} distribute capacity more evenly across services while still accounting for degradation. Notably, the \textit{PL cal. \& cyc.} model shows a stronger emphasis on DCH and baseline discharge in the earlier years, consistent with a preference for lower-SoC operation under calendar ageing. This concentration gradually weakens over time, with service allocation becoming more balanced in later years.

Figs.~\ref{fig:service_year} and \ref{fig:service_ageing} further illustrate the temporal evolution of optimal service allocation. In contrast to ageing-aware strategies, the \textit{No ageing} and \textit{Cycle limit} models exhibit relatively stable service allocations over time, with only minor variations across years. This reflects the absence or weak influence of degradation constraints, resulting in service allocation that is largely insensitive to both time and past battery degradation. By comparison, ageing-aware strategies, particularly the L-based models, show significant changes in service allocation over time. For example, in Year 1, the \textit{L cyc.} and \textit{L cal. \& cyc.} models allocate approximately 66\% and 86\% of capacity to DCH services, respectively, while the \textit{PL cal. \& cyc.} model allocates around 27\%. However, this concentration gradually weakens as the battery ages over time. In the L-based models (\textit{L cyc.} and \textit{L cal. \& cyc.}), the share of DCH decreases substantially as the battery ages, dropping from around 66\% to 25\% in the \textit{L cyc.} model and from about 86\% to 25\% in the \textit{L cal. \& cyc.} model. This represents a reduction of approximately 62\% and 71\%, respectively, respectively, with baseline discharge and other service categories increasing in share, indicating a shift toward more diversified optimal service stacking. In contrast, the piecewise linear models exhibit more moderate changes over time. For example, in the \textit{PL cal. \& cyc.} model, DCH decreases from approximately 27\% in Year 1 to around 17\% in Year 20, while DCL and baseline discharge increase slightly. This reflects a more stable and balanced service allocation throughout the battery lifetime. 

The share of high- and low-frequency services also evolves over time, as shown in Fig.~\ref{fig:service_H_L}. In the L-based models, the share of high-frequency services is initially very high but decreases sharply in the early years, followed by more gradual adjustments thereafter. In contrast, the piecewise linear models exhibit more differentiated behaviour. The \textit{PL cyc.} model remains relatively stable over time, with only minor variations in the high-frequency share. By comparison, the \textit{PL cal. \& cyc.} model shows a noticeable adjustment in the early years, with the high-frequency share decreasing from around 80\% in Year 1 to approximately 58\% in Year 5, before stabilising at later stages. 

These trends indicate that ageing-aware operation affects not only lifetime value, but also the allocation of battery capacity across services over time. Because products differ in delivery duration, response characteristics, and energy-throughput requirements, the preferred service-stacking strategy depends on the operator’s objective, ranging from short-term revenue extraction to lifetime profit maximisation. Overall, the case study shows that optimal service stacking should be treated as a lifetime operational decision, shaped jointly by market revenues, service requirements, degradation costs, and operator objectives.

\section{Conclusion}\label{sec:concl}

Battery storage operators face a coupled lifetime decision problem when stacking multiple services. Participation decisions are linked through shared power and energy capacity, baseline scheduling, and system-operator energy management rules. At the same time, different services impose different technical requirements, activation patterns, and degradation impacts, so their long-term value depends on both market conditions and the evolving health of the battery. Determining the lifetime profit-maximising stacking strategy therefore requires a co-optimisation framework that jointly optimises service participation while accounting for battery ageing and long-term asset value.

We study this problem using Great Britain as a case study, where the new Enduring Auction Capability platform provides a practical setting for analysing coupled multi-service battery operation. Under this platform, Dynamic Containment, Dynamic Moderation, and Dynamic Regulation services are procured within a unified framework, creating new opportunities for revenue stacking while also tightening the coupling between service provision, baseline energy scheduling, and state-of-energy management requirements imposed by the National Energy System Operator. For battery operators, this makes lifetime operating strategy, optimal service stacking, and long-term asset value substantially more complex to assess.

This work develops a lifetime assessment framework for grid-scale battery co-optimisation in this market setting. The framework jointly optimises participation across eight services, while explicitly modelling product-specific technical characteristics and system-operator energy management requirements. Using real GB price and system frequency data, we evaluate battery lifetime performance under a range of operating strategies, from degradation-unaware and rule-based approaches to ageing-aware optimisation with different levels of model fidelity.

Although the analysis is conducted in the GB context, the new Enduring Auction Capability platform demonstrates several advanced market design features, including unified clearing and integrated co-optimisation of services, and therefore provides a useful reference for the future development of frequency response markets in other power systems. In this setting, the proposed framework offers a transferable tool for evaluating how degradation, operational constraints, and market design shape service-stacking decisions across coupled energy and frequency response services, and for identifying strategies that maximise long-term economic value.

The case study highlights three main findings with direct implications for both battery operators and market designers. First, ageing-aware modelling significantly changes the stacked economic value of storage. Ignoring degradation leads to aggressive operation, premature end-of-life, and lower lifetime value, while overly conservative formulations preserve lifetime at the expense of revenue generation. The highest value is achieved when both calendar and cycling ageing are represented with sufficient fidelity. Importantly, calendar ageing exhibits a dual effect: depending on model fidelity, its inclusion can either increase or reduce total lifetime revenue, but it consistently improves the economic value extracted per unit of degradation. This implies that battery operators require sufficiently accurate degradation models not only to balance utilisation and asset preservation, but also to correctly capture the economic role of calendar ageing in long-term operation.

Second, the economically preferred strategy depends on the discount rate. Accounting for both calendar and cycling ageing gives the most complete degradation-aware representation and remains optimal for discount rates up to 15\%. At very high discount rates, the performance gap between the calendar-and-cycle and cycle-only ageing strategies becomes small, with the cycle-only strategy marginally preferred in the numerical results. This occurs because very high discounting reduces the value of preserving battery capacity for future revenues, together with the limited ability of the piecewise-linear approximation to represent nonlinear calendar-ageing behaviour.

Third, the case study determines optimal service-stacking strategies across two electricity arbitrage services and six dynamic frequency response services, linking the service mix to different operating strategies. Across all strategies, high-frequency products are consistently preferred because they create a double revenue opportunity: high-frequency service provision charges the battery and enables subsequent baseline discharge in the energy market. However, the optimal service mix depends on the operator’s economic objective. Operators focused on short-term revenue should adopt a diversified, high-utilisation portfolio in which baseline discharge dominates and capacity is distributed across multiple frequency response services. By contrast, operators seeking to maximise long-term revenue should prioritise high-frequency products, particularly rapid-response Dynamic Containment high-frequency services combined with baseline discharge, while allocating less capacity to more energy-intensive services such as Dynamic Moderation and Dynamic Regulation. Moreover, the optimal service stack is not static. As the battery ages, the degradation rate under the same operating conditions decreases before end-of-life, so the optimal capacity allocation across services also changes over time. In this case, the strategy should focus on Dynamic Containment high-frequency products in early operation, before gradually shifting towards a more balanced service mix as the battery ages.

Future work will focus on extending the framework to incorporate stochastic or chance-constrained optimisation in order to model uncertainty in frequency activation and market conditions, thereby accounting more realistically for forecasting errors and operator risk preferences. In addition, expanding the revenue stack beyond arbitrage and dynamic frequency response to include participation in the balancing mechanism and net imbalance volume chasing strategies would provide a more comprehensive assessment of battery economic value.


\nomenclature{BESS}{Battery Energy Storage System}
\nomenclature{DC}{Dynamic Containment}
\nomenclature{DM}{Dynamic Moderation}
\nomenclature{DR}{Dynamic Regulation}
\nomenclature{DFR}{Dynamic Frequency Response}
\nomenclature{DCH/DCL}{Dynamic Containment in the high-/low-frequency direction}
\nomenclature{DMH/DML}{Dynamic Moderation in the high-/low-frequency direction}
\nomenclature{DRH/DRL}{Dynamic Regulation in the high-/low-frequency direction}
\nomenclature{EA}{Energy Arbitrage}
\nomenclature{EoL}{End of Life}
\nomenclature{GB}{Great Britain}
\nomenclature{IDM}{Intra-day Market}
\nomenclature{NESO}{National Energy System Operator}
\nomenclature{PL}{Piecewise Linear}
\nomenclature{SoC}{State of Charge}
\nomenclature{SoE}{State of Energy}
\nomenclature{SoS2}{Special Ordered Sets of Type 2}

\renewcommand{\nomname}{List of Abbreviations}
\printnomenclature

\appendix
\section*{Nomenclature}
\subsection*{Sets}
\begin{description}[leftmargin=!, labelwidth=3.3cm]
    \setlength{\itemsep}{0pt}
    \setlength{\parskip}{0pt}
   \item[$\mathcal{H}^{\mathrm{15m}}_{\Delta t}(t)$] 
    Index set of timesteps at resolution $\Delta t$ contained in the 15-minute time window starting at $t \in \mathcal{T}_{\mathrm{15m}}$
    \item[$\mathcal{H}^{\mathrm{4h}}_{\Delta t}(t)$] 
    Index set of timesteps at resolution $\Delta t$ contained in the 4-hour time window starting at $t \in \mathcal{T}_{\mathrm{4h}}$
    \item[$\mathcal{H}^{\mathrm{daily}}_{\Delta t}(t)$] 
    Index set of timesteps at resolution $\Delta t$ contained in the daily time window starting at $t \in \mathcal{T}_{\mathrm{daily}}$
    \item[$\mathcal{H}^{\mathrm{EFA}}_{\Delta t}(t)$] 
    Index set of timesteps at resolution $\Delta t$ contained in the Electricity Forward Agreement block starting at $t \in \mathcal{T}_{\mathrm{EFA}}$
    \item[$\mathcal{H}^{\mathrm{EFA}}_{\mathrm{SP}}(t)$] 
    Index set of settlement periods covered by the Electricity Forward Agreement block starting at $t \in \mathcal{T}_{\mathrm{EFA}}$
    \item[$\mathcal{H}^{\mathrm{SP}}_{\Delta t}(t)$] 
    Index set of timesteps at resolution $\Delta t$ contained in the settlement period starting at $t \in \mathcal{T}_{\mathrm{SP}}$
    \item[$\mathcal{J}$] Index set of breakpoints for the piecewise-linear calendar ageing model
    \item[$\mathcal{K}$] Index set of breakpoints for the piecewise-linear cycle ageing model
    \item[$\mathcal{S}$] Index set of dynamic frequency response services, $\mathcal{S}=\{\mathrm{DCH},\mathrm{DCL},\mathrm{DMH},\mathrm{DML},\mathrm{DRH},\\ \mathrm{DRL}\}$
    \item[$\mathcal{T}$] Index set of timesteps at resolution $\Delta t$ for the current optimisation horizon
    \item[$\mathcal{T}_{\mathrm{15m}}$] Index set of 15-minute time intervals over the current optimisation horizon
    \item[$\mathcal{T}_{\mathrm{4h}}$] Index set of 4-hour time intervals over the current optimisation horizon
    \item[$\mathcal{T}_{\mathrm{daily}}$] Index set of daily time intervals over the current optimisation horizon
    \item[$\mathcal{T}_{\mathrm{EFA}}$] Index set of Electricity Forward Agreement blocks over the current optimisation horizon
    \item[$\mathcal{T}_{\mathrm{EFA}}^{\mathrm{SP1}}$] Index set of settlement periods corresponding to the first settlement period of each Electricity Forward Agreement block over the current optimisation horizon
    \item[$\mathcal{T}_{\mathrm{EFA}}^{\mathrm{SP8}}$] Index set of settlement periods corresponding to the last settlement period of each Electricity Forward Agreement block over the current optimisation horizon
    \item[$\mathcal{T}_{\mathrm{EFA}}^{\mathrm{SP2\text{--}8}}$] Index set of settlement periods excluding the first settlement period of each Electricity Forward Agreement block over the current optimisation horizon
    \item[$\mathcal{T}_{\mathrm{SP}}$] Index set of settlement periods over the current optimisation horizon
\end{description}

\subsection*{Parameters}
\begin{description}[leftmargin=!, labelwidth=3.2cm]
    \setlength{\itemsep}{0pt}
    \setlength{\parskip}{0pt}
        \item[$a_{\mathrm{cal}}$] Linear calendar-ageing coefficient (slope) 
        \item[$a_{\mathrm{cyc}}$] Linear cycle-ageing coefficient (slope)

        \item[$a^{\mathrm{cyc,ch/dis}}$] Linear cycle-ageing coefficient for charging/discharging (slope)
        \item[$b^{\mathrm{cyc,ch/dis}}$] Linear cycle-ageing coefficient for charging/discharging (intercept)
        \item[$\alpha_t^{s}$] DFR activation factor for service $s$ in timestep $t$, based on actual system frequency
        \item[$b_{\mathrm{cal}}$] Linear calendar-ageing coefficient (intercept)
        \item[$b_{\mathrm{cyc}}$] Linear cycle-ageing coefficient (intercept)
        \item[$\Delta t$] Duration of a single optimisation timestep
        \item[$E^{\mathrm{batt}}$] Energy capacity of the battery
        \item[$\eta$] Charge/discharge efficiency of the BESS assumed for the optimisation model
        \item[$\Lambda^{\mathrm{SP}}$] Duration of a settlement period, fixed at 0.5 hours
        \item[$\Lambda^{\mathrm{EFA}}$] Duration of an electricity forward agreement block, fixed at 4 hours
        \item[$n^{\mathrm{cyc,max}}$] Max number of full cycle equivalent per day
        \item[$P^{\max}$] Maximum charge and discharge power of the BESS
        \item[$\pi^{\mathrm{IDM}}_t$] Intra-day electricity price from EPEX SPOT in timestep $t$
        \item[$\pi^{\mathrm{loss}}$] Ageing cost per unit of capacity loss for the optimisation model 
        \item[$\pi^{s}_t$] Clearing price of frequency service $s$ from NESO in timestep $t$ 
        \item[$\rho$] Scaling factor linking frequency response energy to baseline charging and discharging actions
        \item[$\mathrm{SoH}^{\mathrm{EoL}}$] End-of-life SoH threshold, fixed at 80\%
        \item[$\mathrm{SoC}^{\mathrm{ini}}$] Initial SoC at the beginning of the optimisation horizon 
        \item[$X_j^{\mathrm{cal}/\mathrm{cyc}}$] Breakpoint values of the calendar/cycle ageing functions at point $j$
        \item[$Z_j^{\mathrm{cal}/\mathrm{cyc}}$] Degradation values of the calendar/cycle ageing functions at point $j$
\end{description}

\subsection*{Variables}
\begin{description}[leftmargin=!, labelwidth=3.2cm]
    \setlength{\itemsep}{0pt}
    \setlength{\parskip}{0pt}
  \item[$\mathrm{ABS}_t^{\mathrm{H}},\, \mathrm{ABS}_t^{\mathrm{L}}$] 
Energy absorbed for high- and low-frequency response at the start of settlement period $t$
\item[$\mathrm{adj0}_t^{\mathrm{H}},\, \mathrm{adj0}_t^{\mathrm{L}}$] 
Energy recovery adjustment volume scheduled for high- and low-frequency response at the start of settlement period $t$
\item[$\mathrm{adj4}_t^{\mathrm{H}},\, \mathrm{adj4}_t^{\mathrm{L}}$] 
Energy recovery adjustment volume implemented for high- and low-frequency response at the start of settlement period $t$

\item[$b_t^{\mathrm{H}},\, b_t^{\mathrm{L}}$] 
Binary indicators for provision of high- and low-frequency response services in timestep $t$

\item[$C^{\mathrm{ageing}}$] 
Total battery degradation cost over the optimisation horizon

\item[$\mathrm{CREV}_t^{\mathrm{H}},\, \mathrm{CREV}_t^{\mathrm{L}}$] 
Contracted response energy volume for high- and low-frequency services over the Electricity Forward Agreement block starting at timestep $t$

\item[$e^{\mathrm{ch}}_t,\,e^{\mathrm{dis}}_t$] 
Energy throughput in the charging/discharging direction during the 4-hour window in timestep $t$

\item[$\mathrm{ER}_t^{\mathrm{H}},\, \mathrm{ER}_t^{\mathrm{L}}$] 
Energy recovery volume for the high- and low-frequency direction over the Electricity Forward Agreement block starting at timestep $t$

\item[$\mathrm{FRE}_t^{\mathrm{H}},\, \mathrm{FRE}_t^{\mathrm{L}}$] 
Frequency response energy delivered for high- and low-frequency services during settlement period $t$

\item[$\mathrm{left}_t^{\mathrm{H}},\, \mathrm{left}_t^{\mathrm{L}}$] 
Accumulated unrecovered energy from high- and low-frequency response at the start of settlement period $t$

\item[$\mathrm{MG}_t^{\mathrm{H}},\, \mathrm{MG}_t^{\mathrm{L}}$] 
State-of-energy margin for high- and low-frequency response at the start of settlement period $t$

\item[$\mathrm{MSER}_t^{\mathrm{H}},\, \mathrm{MSER}_t^{\mathrm{L}}$] 
Minimum state-of-energy required to deliver high- and low-frequency response at the start of settlement period $t$

\item[$n_t^{\mathrm{cyc}}$] 
Number of equivalent full charge--discharge cycles accumulated over day $t$

\item[$p_t^{\mathrm{base,ch}},\,p_t^{\mathrm{base,dis}}$] 
Baseline charging/discharging power in timestep $t$

\item[$p_t^{s}$] 
Contracted power for DFR service $s \in \mathcal{S}$ in timestep $t$

\item[$p_t^{\mathrm{tot}}$] 
Total power exchanged with the grid in timestep $t$

\item[$p_t^{\mathrm{tot,ch}},\,p_t^{\mathrm{tot,dis}}$] 
Total charging/discharging power in timestep $t$

\item[$\Pi^{\mathrm{DFR}}$] 
Revenue from dynamic frequency response services

\item[$\Pi^{\mathrm{E}}$] 
Revenue from the electricity market

\item[$q_t^{\mathrm{loss,cal}},\,q_t^{\mathrm{loss,cyc}}$] 
Capacity loss due to calendar/cycle ageing in timestep $t$

\item[$q_t^{\mathrm{loss,cyc,ch}},\, q_t^{\mathrm{loss,cyc,dis}}$] 
Capacity loss from charging/discharging-related cycle ageing in timestep $t$

\item[$r_t^{\mathrm{H}},\, r_t^{\mathrm{L}}$] 
Reserved power for energy recovery for high- and low-frequency services in timestep $t$

\item[$\mathrm{RER}_t^{\mathrm{H}},\, \mathrm{RER}_t^{\mathrm{L}}$] 
Required energy recovery in the high- and low-frequency direction during settlement period $t$, representing the maximum recoverable energy in period $t$ subject to the recovery limit

\item[$\mathrm{SoC}_t$] 
State of charge of the battery in timestep $t$

\item[$\mathrm{SoE}_t$] 
State of energy in timestep $t$

\item[$\overline{\mathrm{SoC}}_t$] 
Average state of charge over the 15-minute window in timestep $t$

\item[$z_t^{\mathrm{base,ch}},\,z_t^{\mathrm{base,dis}}$] 
Binary variable indicating whether baseline charging/discharging is active in timestep $t$

\item[$z_t^{\mathrm{tot,ch}},\, z_t^{\mathrm{tot,dis}}$] 
Binary variable indicating whether total charging/discharging is active in timestep $t$

\item[$\lambda^{\mathrm{cal}}_{t,j}$] 
SOS-type~2 variables for linearisation of calendar ageing

\item[$\lambda^{\mathrm{cyc,ch}}_{t,k},\, \lambda^{\mathrm{cyc,dis}}_{t,k}$] 
SOS-type~2 variables for linearisation of cycle ageing in the charging/discharging direction

\end{description}

\section{Illustrative Examples of NESO SoE Management Rules} \label{app:neso_soe}
NESO imposes SoE management requirements on DFR service providers to ensure sustained delivery capability over the contracted service period \cite{NESO_SOE_monitor,NGESO_SOE_ppt}. Central to this framework is the energy recovery mechanism, which governs SoE recovery following DFR service activation. The energy recovery requirement is defined as 20\% of the contracted response energy volume in MWh and represents the maximum energy volume that may need to be restored per settlement period following a shortfall relative to the threshold. When frequency service activation causes the SoE to fall below the minimum state-of-energy requirement level, providers must submit an operational baseline adjustment in the electricity market at the earliest permissible opportunity. The corresponding recovery volume, referred to as the energy recovery adjustment volume, adjusts the minimum SoE requirement and therefore affects the SoE trajectory. It is determined by the realised frequency response energy and the prevailing SoE margin relative to the contracted response energy volume threshold \cite{NESO_new_explain,UKbattery_SoE_mana}.

To ensure that post-event recovery remains feasible in real time following service delivery, NESO imposes power reserve requirements at the \emph{day-ahead} auction stage. For all three DFR services, units submitting bi-directional bids are required to reserve a specified proportion of their power capacity for energy recovery, whereas this requirement does not apply to units offering single-directional services. Specifically, providers with bi-directional contracts must reserve additional power capacity equivalent to 10\%, 20\%, and 40\% of the contracted power in MW for DC, DM, and DR services, respectively\footnote{The energy recovery requirement is defined in energy terms (MWh) and is mapped to a per settlement period power reserve by converting the energy volume to its 30-minute power equivalent (multiplication by 2) and scaling by the service-specific delivery duration. For example, a unit with a 40 MW DC contract implies a contracted response energy volume of 10 MWh. Applying the 20\% energy recovery requirement yields a reserve energy of 2 MWh per settlement period, corresponding to a charging power of 4 MW over a 30-minute interval. This translates into a reserve requirement equal to 10\% of the contracted DC capacity at the bidding stage.}. This power reserve requirement ensures that bi-directional DFR service providers retain sufficient operational flexibility to manage their SoE across all contracted service periods. By maintaining adequate headroom and footroom, the battery can deliver the full contracted response energy volume for both high- and low-frequency services. It can also recover energy in time to support continued service provision within the current EFA block and meet the entry requirements of subsequent blocks.

\begin{figure}[t] 
  \centering
  \includegraphics[width=\columnwidth]{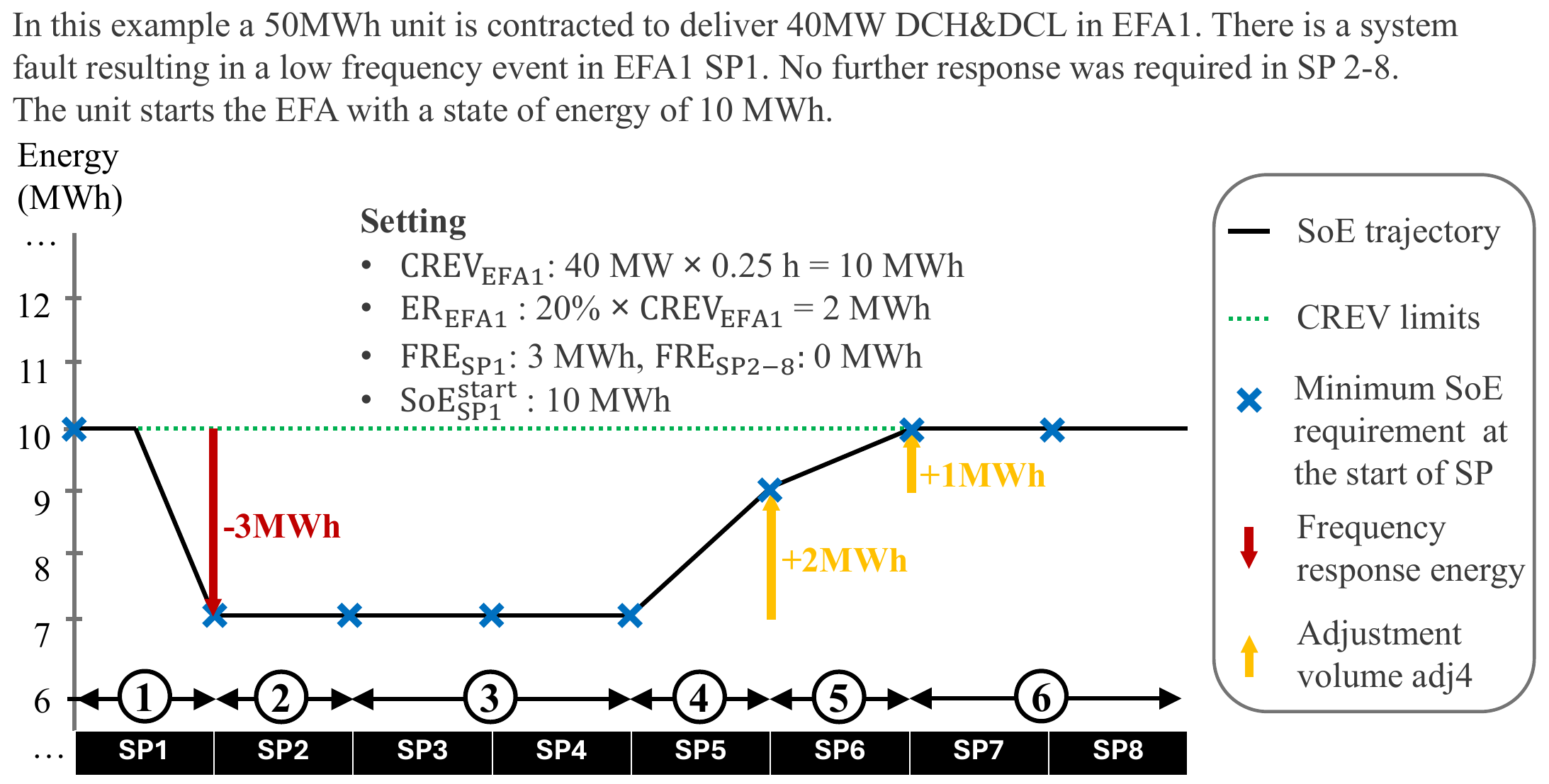}
  \caption{NESO SoE management example 1: Initial SoE = 10 MWh. DCH: dynamic containment high-frequency services; DCH: dynamic containment low-frequency services; EFA: electricity forward agreement; SP: settlement period; SoE: state-of-energy; CREV: contracted response energy volume; ER: energy recovery; FRE: frequency response energy.}\label{fig:neso_case1}
\end{figure}

\begin{figure}[t] 
  \centering
  \includegraphics[width=\columnwidth]{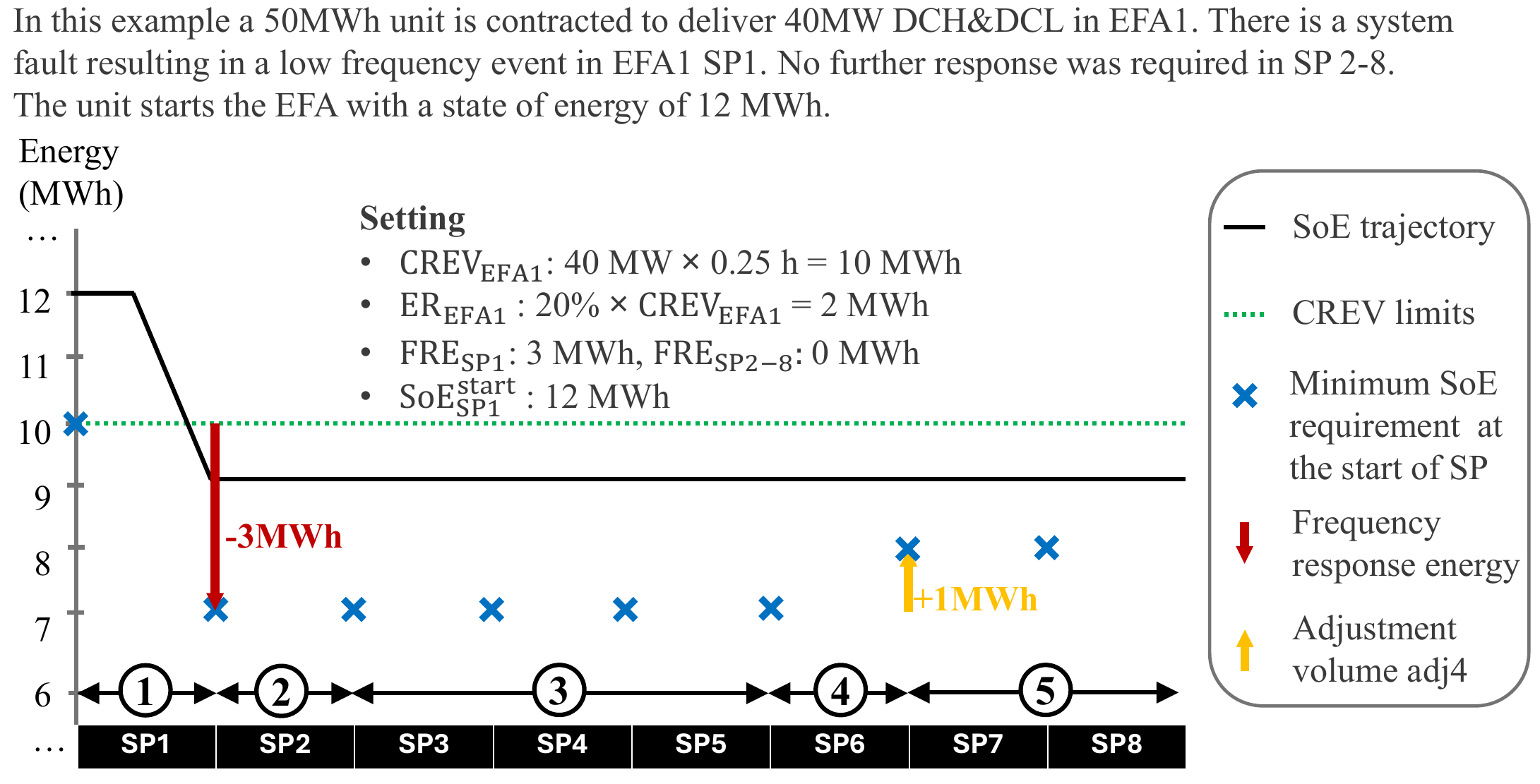}
  \caption{NESO SoE management example 2: Initial SoE = 12 MWh. DCH: dynamic containment high-frequency services; DCH: dynamic containment low-frequency services; EFA: electricity forward agreement; SP: settlement period; SoE: state-of-energy; CREV: contracted response energy volume; ER: energy recovery; FRE: frequency response energy.}\label{fig:neso_case2}
\end{figure}

Figs.~\ref{fig:neso_case1} and~\ref{fig:neso_case2} present two illustrative examples for the SoE management rules, intended to be read alongside Section~\ref{sec:mpc:opt:soe_mana} and Table~\ref{tab:soe_management_abbreviations} to aid understanding of how the constraints operate in practice. In both examples, a 50 MWh unit is contracted to deliver 40 MW each of DCH and DCL services within electricity forward agreement EFA1. A system fault results in a low-frequency event in the first settlement period, triggering contract activation and a delivered frequency response energy of 3 MWh, with no further response required in settlement periods 2--8. The initial SoE at the start of the EFA1 is set to 10 MWh in Case 1 and 12 MWh in Case 2, while all other settings remain the same. In both examples, the contracted response energy volume is 10 MWh for each direction, calculated as 40 MW $\times$ 0.25 h = 10 MWh, where 0.25 h (15 minutes) represents the delivery duration of DC services. This corresponds to 2 MWh of energy recovery capacity (20\% of contracted response energy volume) derived from reserved power at the bidding stage. The minimum SoE requirement at the start of settlement period 1 equals the contracted response energy volume ($10\,\text{MWh}$) and evolves cumulatively across settlement periods. The minimum SoE requirement is calculated at the start of each settlement period and shall: (i) equal the contracted response energy volume at the start of the contracted service period; (ii) decrease with service activation, i.e., frequency response energy; and (iii) increase with the implemented energy recovery adjustment volume, thereby moving back towards the contracted response energy volume (see Constraints~\eqref{eq:neso_soe_mana:MSER_START}--\eqref{eq:neso_soe_mana:MSER_cal_H}). The energy recovery adjustment volume is applied to the minimum SoE requirement at the start of each settlement period, rather than directly to the response unit’s SoE. The unit must maintain its SoE no less than the minimum SoE requirement at the start of every settlement period of contracted service period, but has more flexibility with regard to how it manages its SoE to stay within the allowable range.

\begin{table}[t]
\centering
\footnotesize
\caption{Values of SoE management variables across service periods under NESO SoE management (Example 1)}
\begin{tabular}{p{8cm}cccccccc}
\toprule
\midrule
 & SP1 & SP2 & SP3 & SP4 & SP5 & SP6 & SP7 & SP8 \\
\midrule
SoE at the start of the settlement period ($\mathrm{SoE}^{\mathrm{start}}_{\mathrm{sp}}$) & 10 & 7 & 7 & 7 & 7 & 9 & 10 & 10 \\
Leftover at the start of the settlement period ($\mathrm{left}_{\mathrm{sp}}$) & 0 & 1 & 0 & 0 & 0 & 0 & 0 & 0 \\
Minimum SoE requirement at the start of the settlement period ($\mathrm{MSER}_{\mathrm{sp}}$) & 10 & 7 & 7 & 7 & 7 & 9 & 10 & 10 \\
Minimum SoE requirement at the end of the settlement period ($\mathrm{MSER}_{\mathrm{sp+1}}$) & 7 & 7 & 7 & 7 & 9 & 10 & 10 & / \\
Frequency response energy during the settlement period ($\mathrm{FRE}_{\mathrm{sp}}$) & 3 & 0 & 0 & 0 & 0 & 0 & 0 & 0 \\
SoE margin at the start of the settlement period ($\mathrm{MG}_{\mathrm{sp}}$) & 0 & 0 & 0 & 0 & 0 & 0 & 0 & 0 \\
Required energy recovery of the settlement period ($\mathrm{RER}_{\mathrm{sp}}$) & 2 & 1 & 0 & 0 & 0 & 0 & 0 & 0 \\
Energy absorbed of the settlement period ($\mathrm{ABS}_{\mathrm{sp}}$) & 0 & 0 & 0 & 0 & 0 & 0 & 0 & 0 \\
Energy recovery adjustment volume scheduled of the settlement period ($\mathrm{adj0}_{\mathrm{sp}}$) & 2 & 1 & 0 & 0 & 0 & 0 & 0 & 0 \\
Energy recovery adjustment volume implemented of the settlement period ($\mathrm{adj4}_{\mathrm{sp}}$) & 0 & 0 & 0 & 0 & 2 & 1 & 0 & 0 \\
\midrule
\bottomrule
\end{tabular}
\label{tab:neso_case1}
\end{table}

In the first example, as illustrated in Fig.~\ref{fig:neso_case1} and reported in Table~\ref{tab:neso_case1}, the initial SoE at the start of the EFA1 is set to 10 MWh, equal to the minimum SoE requirement of 10 MWh. The settlement period 1 margin is therefore zero. In stage \circled{1} (settlement period 1), the frequency event causes a DCL service activation of 3 MWh. At the end of settlement period 1 (i.e., at the start of settlement period 2), the minimum SoE requirement drops to 7 MWh (10 MWh + 0 MWh $-$ 3 MWh) due to service activation and the absence of energy recovery adjustment volume implemented for settlement period 1. The required energy recovery volume is calculated as the lower of: (i) the frequency response energy plus any leftover, equal to 3 MWh; and (ii) energy recovery capability per settlement period (20\% of contracted response energy volume), equal to 2 MWh. The required energy recovery for settlement period 1 is therefore 2 MWh. The absorption is defined as the minimum of required energy recovery and SoE margin, which is zero. The energy recovery adjustment volume implemented in settlement period 5, after the 4-settlement-period implementation delay, is calculated as the difference between the required energy recovery and energy absorption, yielding 2 MWh. This implies that only 2 MWh of the 3 MWh service delivery is scheduled for adjustment in settlement period 5, as the energy recovery capability limits the adjustment that can be made per settlement period. Consequently, the leftover 1 MWh (3 MWh - 2 MWh) must be recovered in subsequent settlement periods. In stage \circled{2} (settlement period 2), the SoE is required to be at least the minimum SoE requirement of 7 MWh at the start of settlement period 2. No service is activated during this settlement period; frequency response energy is zero in settlement period 2 and remains zero in subsequent settlement periods. The margin is also zero. As no energy recovery adjustment volume is implemented and no frequency response energy is delivered in the current settlement period, the minimum SoE requirement remains unchanged at 7 MWh at the end of settlement period 2. The required energy recovery, however, equals 1 MWh, reflecting the 1 MWh leftover to be recovered. Accordingly, settlement period 2 sets an energy recovery adjustment volume of 1 MWh, which will be implemented in settlement period 2+4 (settlement period 6). This energy recovery adjustment volume fully accounts for the remaining 1 MWh of energy to be recovered. In stage \circled{3} (settlement periods 3–4), the minimum SoE requirement remains unchanged at 7 MWh for the same reason as in stage \circled{2}, and no energy recovery adjustment volume is determined or implemented. In stage \circled{4} (settlement period 5), the energy recovery adjustment volume of 2 MWh determined in stage \circled{1} (settlement period 1) is implemented, and the minimum SoE requirement therefore increases by 2 MWh during settlement period 5. This is achieved by submitting appropriate orders in the intraday market. Similarly, in stage \circled{5} (settlement period 6), the energy recovery adjustment volume of 1 MWh determined in stage \circled{2} (settlement period 2) is implemented, resulting in a further 1 MWh increase in the minimum SoE requirement. As a result, the minimum SoE requirement moves back towards the contracted response energy volume. Finally, in stage \circled{6} (settlement periods 7–8), the minimum SoE requirement remains at the contracted response energy volume value, and the minimum SoE requirement at the end of settlement period 8 equals the contracted response energy volume for the subsequent electricity forward agreement block.

\begin{table}[t]
\centering
\footnotesize
\caption{Values of SoE management variables across service periods under NESO SoE management (Example 2)}
\begin{tabular}{p{8cm}cccccccc}
\toprule
\midrule
 & SP1 & SP2 & SP3 & SP4 & SP5 & SP6 & SP7 & SP8 \\
\midrule
SoE at the start of the settlement period ($\mathrm{SoE}^{\mathrm{start}}_{\mathrm{sp}}$) & 12 & 9 & 9 & 9 & 9 & 9 & 9 & 9 \\
Leftover at the start of the settlement period ($\mathrm{left}_{\mathrm{sp}}$) & 0 & 1 & 0 & 0 & 0 & 0 & 0 & 0 \\
Minimum SoE requirement at the start of the settlement period ($\mathrm{MSER}_{\mathrm{sp}}$) & 10 & 7 & 7 & 7 & 7 & 7 & 8 & 8 \\
Minimum SoE requirement at the end of the settlement period ($\mathrm{MSER}_{\mathrm{sp+1}}$) & 7 & 7 & 7 & 7 & 7 & 8 & 8 & / \\
Frequency response energy during the settlement period ($\mathrm{FRE}_{\mathrm{sp}}$) & 3 & 0 & 0 & 0 & 0 & 0 & 0 & 0 \\
SoE margin at the start of the settlement period ($\mathrm{MG}_{\mathrm{sp}}$) & 2 & 0 & 0 & 0 & 0 & 0 & 0 & 0 \\
Required energy recovery of the settlement period ($\mathrm{RER}_{\mathrm{sp}}$) & 2 & 1 & 0 & 0 & 0 & 0 & 0 & 0 \\
Energy absorbed of the settlement period ($\mathrm{ABS}_{\mathrm{sp}}$) & 2 & 0 & 0 & 0 & 0 & 0 & 0 & 0 \\
Energy recovery adjustment volume scheduled of the settlement period ($\mathrm{adj0}_{\mathrm{sp}}$) & 0 & 1 & 0 & 0 & 0 & 0 & 0 & 0 \\
Energy recovery adjustment volume implemented of the settlement period ($\mathrm{adj4}_{\mathrm{sp}}$) & 0 & 0 & 0 & 0 & 0 & 1 & 0 & 0 \\
\midrule
\bottomrule
\end{tabular}
\label{tab:neso_case2}
\end{table}

In the second example, as illustrated in Fig.~\ref{fig:neso_case2} and reported in Table~\ref{tab:neso_case2}, the initial SoE at the start of the electricity forward agreement is 12 MWh, exceeding the minimum SoE requirement of 10 MWh and resulting in a settlement period 1 margin of 2 MWh. As in Case 1, a 3 MWh DCL activation in stage \circled{1} (settlement period 1) reduces the minimum SoE requirement to 7 MWh at the start of settlement period 2, and the required energy recovery equals 2 MWh (the minimum of 2 MWh and 3 MWh). Unlike Case 1, the absorption equals the margin (2 MWh), leading to an energy recovery adjustment volume of 0 MWh for settlement period 1+4 (settlement period 5). Consequently, no energy adjustment is required in settlement period 5, as the additional initial headroom implies that the SoE is expected to remain at least at the minimum SoE requirement level at the end of settlement period 5. This flexibility associated with a higher starting margin in SoE management has been recently introduced to allow units to recover less energy later in the contracted service period \cite{NESO_new_explain}. The remaining 1 MWh (3 MWh - 2 MWh) is recovered in subsequent settlement periods. In stage \circled{2} (settlement period 2), the required energy recovery equals 1 MWh, reflecting the 1 MWh leftover to be recovered despite the 2 MWh SoE headroom at the start of settlement period 1. However, the energy margin remains zero from settlement period 2 onwards, and settlement period 2 sets an energy recovery adjustment volume of 1 MWh, which will be implemented in settlement period 2+4 (settlement period 6). In stage \circled{3} (settlement periods 3–5), the minimum SoE requirement remains unchanged at 7 MWh and no adjustment is implemented, in contrast to Case 1 where an energy recovery adjustment volume of 2 MWh is implemented in settlement period 5. In stage \circled{4} (settlement period 6), the energy recovery adjustment volume of 1 MWh determined in stage \circled{2} (settlement period 2) is implemented, resulting in a 1 MWh increase in the minimum SoE requirement. Notice that the minimum SoE requirement does not return to the contracted response energy volume within the current contracted service period because part of the required recovery is absorbed by the initial SoE headroom, resulting in a reduced energy recovery adjustment volume being applied to the minimum SoE requirement. Importantly, the revised rules do not relax the minimum SoE requirement; instead, they provide flexibility in how units manage their SoE to remain no less than the evolving threshold. Finally, in stage \circled{5} (settlement periods 7–8), the minimum SoE requirement remains at 8 MWh, and the minimum SoE requirement at the end of settlement period 8 equals the contracted response energy volume for the subsequent electricity forward agreement block.

\section{Key parameters for the model predictive control framework}\label{sec:para_setting}
Key parameters for the model predictive control framework are summarised in Table~\ref{tab:model_params}.
\begin{table}[htbp]
\centering
\caption{Model Parameters and Assumptions.}
\setlength{\tabcolsep}{3.5pt} 
\begin{tabularx}{\textwidth}{l X X p{1.3cm} p{2cm}}
\toprule
\midrule
\textbf{Category} & \textbf{Parameter} & \textbf{Value/Metric} & \textbf{Unit} & \textbf{Source} \\
\midrule

\multirow{6}{*}{\makecell[tl]{Optimisation\\Configuration}} 
 & $\Delta t$ & 60 & s & - \\
 & $\mathcal{T}$ & 2 & days & - \\
 & MIP Gap & 1 & \% & - \\
 & Time limit per optimisation & 1,800 & s & - \\
\midrule

\multirow{5}{*}{\makecell[tl]{Simulation\\Configuration}}
 & AC/DC converter & - & - & \cite{NOTTON_converter} \\
 & Cell Type & SonyLFP (Sony US26650FTC1) & - & \cite{SimSES_paper,calen_ageing_LFP,cyc_ageing_LFP} \\
 & $P^{\max}$ & 5 & MW & - \\
 & $E^{\text{batt}}$ & 5 & MWh & - \\
 & Initial SoE & 50 & \% & - \\
\midrule

\multirow{4}{*}{\makecell[tl]{Optimisation\\Parameter}}
 & $\eta$ & 0.9 & - & \cite{MPC_battery} \\
 & $\delta$ & 0.5 & - & - \\
 & $n_{\text{cyc,max}}$ & 2 & - & \cite{UKbattery_SoE_mana} \\
 & $\pi^{\mathrm{loss}}$ & 150,000 & \pounds/MWh & \cite{WANKMULLER201756} \\
 & $\mathrm{SoH}^{\text{EoL}}$ & 80 & $\%$ & \cite{MPC_battery} \\

 \midrule
\bottomrule
\end{tabularx}
\label{tab:model_params}
\end{table}

\clearpage
\bibliographystyle{elsarticle-num}
\bibliography{reference}

\end{document}